\documentclass[11pt,a4paper]{article}

\usepackage[T1]{fontenc}
\usepackage[margin=27mm]{geometry}
\usepackage{amsmath,amssymb,amsthm}
\usepackage{abstract}
\usepackage[hidelinks]{hyperref}
\allowdisplaybreaks[2]
\newtheorem{theorem}{Theorem}
\newtheorem{lemma}[theorem]{Lemma}
\newtheorem{proposition}[theorem]{Proposition}
\newtheorem{conjecture}{Conjecture}
\newcommand{\BGW}{\mathrm{BGW}}
\newcommand{\dd}{\mathbf d}
\newcommand{\Rmnum}[1]{\uppercase\expandafter{\romannumeral #1}}  
		 
\newcommand{\Addresses}{{
  \bigskip
  \footnotesize

  Xuanyu Huang, \textsc{Mathematical Science Research Center, Chongqing University of Technology, Chongqing, 400054, China}\par\nopagebreak
  \textit{E-mail address}: \texttt{Hxuanyu98@gmail.com}

  \medskip

  Kefeng Liu, \textsc{Mathematical Science Research Center, Chongqing University of Technology, Chongqing, 400054, China}\par\nopagebreak
  \textit{E-mail address}: \texttt{kefengliu@cqut.edu.cn}

  \medskip

  Hao Xu, \textsc{Mathematical Science Research Center, Chongqing University of Technology, Chongqing, 400054, China}\par\nopagebreak
  \textit{E-mail address}: \texttt{haoxu@cqut.edu.cn}
}}

\title{The super Weil-Petersson asymptotic expansion \\ and the coefficients}
\author{Xuanyu Huang \and Kefeng Liu \and Hao Xu}
\date{}

\begin{document}
\maketitle

\begin{abstract}
We determine the universal constant in the large-genus asymptotic expansion
of super Weil-Petersson volumes. 
The argument adopts the uniform-transport strategy
for shifted theta intersection numbers. The
uniform large-genus theorem for normalized BGW
numbers of Guo-Norbury-Yang-Zagier determines the
amplitude at a simpler endpoint. The result is \(C=\pi^{-9/2}\).
We give a proof of the Griguolo-Papalini-Russo-Seminara conjecture for
fixed boundary lengths and fixed descendant degrees. The boundary
formula is uniform when the sum of the lengths is \(o(\sqrt g)\),
and we exhibit an obstruction at the \(\sqrt g\) scale. We also 
compute the first four inverse-genus coefficients by two exact
recursive methods and derive general polynomial structure results. 
Finally, we show that, for each fixed genus, the super Weil–Petersson volume divided by \(n!\) is a polynomial in the number \(n\) of marked points. 
\end{abstract}

\tableofcontents

\section{Introduction}

Super Weil-Petersson volumes connect the geometry of moduli spaces of
super Riemann surfaces with intersection theory on moduli spaces of
stable curves and with JT supergravity. Stanford and Witten \cite{stanford2020jt} derived a recursion for these volumes as part of
their study of JT gravity and random matrix ensembles. See also \cite{AN,HPZ,kimura2020jt}. Norbury
\cite{norbury2023new,norbury2020enumerative} constructed the cohomology
classes that express the super volumes as intersection numbers and
gave an algebro-geometric proof of the recursion. In the normalization
used here, the volume with geodesic boundary lengths
\(L_1,\ldots,L_n\geq0\) is
\begin{equation}\label{VTheta}
 V^\Theta_{g,n}(L_1,\ldots,L_n)
 =\int_{\overline{\mathcal M}_{g,n}}\Theta_{g,n}
   \exp\left(2\pi^2\kappa_1+\frac12\sum_{i=1}^nL_i^2\psi_i\right).
\end{equation}
Our purpose is to determine the universal constant in the large-genus
expansion of these volumes, prove the resulting asymptotic formulae
for boundary volumes and super intersection numbers, and study their
correction coefficients.

Here \(\overline{\mathcal M}_{g,n}\) is the moduli space of stable
\(n\)-pointed complex curves of genus \(g\), with \(2g-2+n>0\).
Let
\[
 \pi:\overline{\mathcal M}_{g,n+1}\longrightarrow
       \overline{\mathcal M}_{g,n}
\]
be the morphism forgetting the last marking, let \(\sigma_i\) be its
canonical sections, and let \(D_i\) be the corresponding divisors.
If \(\omega_\pi\) denotes the relative dualizing sheaf, the
tautological classes are \cite{arbarello1996combinatorial}
\[
 \psi_i=c_1(\sigma_i^*\omega_\pi),\qquad
 \kappa_j=\pi_*\left(c_1\left(\omega_\pi\left(\sum_iD_i\right)\right)^{j+1}\right).
\]
Denote the gluing morphisms by
\[
\begin{aligned}
 \phi_{\mathrm{irr}}:
 &\overline{\mathcal M}_{g-1,n+2}\longrightarrow
   \overline{\mathcal M}_{g,n},\\
 \phi_{h,I}:
 &\overline{\mathcal M}_{h,|I|+1}\times
   \overline{\mathcal M}_{g-h,|J|+1}\longrightarrow
   \overline{\mathcal M}_{g,n},
 \qquad I\sqcup J=\{1,\ldots,n\},
\end{aligned}
\]
where both factors of a separating boundary stratum are stable.

Norbury introduced classes
\(\Theta_{g,n}\in H^{4g-4+2n}(\overline{\mathcal M}_{g,n},\mathbb Q)\)
, which satisfy the boundary and forgetful-map identities
\[
\begin{aligned}
 \phi_{\mathrm{irr}}^*\Theta_{g,n}&=\Theta_{g-1,n+2},\\
 \phi_{h,I}^*\Theta_{g,n}
 &=\pi_1^*\Theta_{h,|I|+1}\cdot\pi_2^*\Theta_{g-h,|J|+1},\\
 \Theta_{g,n+1}&=\psi_{n+1}\pi^*\Theta_{g,n},
 \qquad \Theta_{1,1}=3\psi_1,
\end{aligned}
\]
where \(\pi_1,\pi_2\) are the projections from the product in the
separating gluing map. In genus zero, the complex degree \(n-2\) of
\(\Theta_{0,n}\) exceeds \(\dim_{\mathbb C}\overline{\mathcal M}_{0,n}=n-3\),
so \(\Theta_{0,n}=0\). In particular, all genus-zero theta
intersection numbers vanish.

For a nonnegative degree vector \(\mathbf d=(d_1,\ldots,d_n)\) with
\(|\mathbf d|=\sum_i d_i\leq g-1\), set \(d_0=g-1-|\mathbf d|\)
and define the normalized super intersection number
\begin{equation}\label{superint}
 [\tau_{d_1}\cdots\tau_{d_n}]^\Theta_g
 :=\frac{2^{2|\mathbf d|}(2\pi^2)^{d_0}}{d_0!}
   \prod_{i=1}^n(2d_i+1)!!
   \int_{\overline{\mathcal M}_{g,n}}
   \Theta_{g,n}\prod_{i=1}^n\psi_i^{d_i}\kappa_1^{d_0}.
\end{equation}
The dimension constraint follows from
\(\deg_{\mathbb C}\Theta_{g,n}=2g-2+n\). Expanding
\eqref{VTheta} gives
\begin{equation}\label{VTheta2}
 V^\Theta_{g,n}(2L_1,\ldots,2L_n)
 =\sum_{|\mathbf d|\leq g-1}
  [\tau_{d_1}\cdots\tau_{d_n}]^\Theta_g
  \prod_{i=1}^n\frac{L_i^{2d_i}}{(2d_i+1)!}.
\end{equation}
We write \(V^\Theta_{g,n}=V^\Theta_{g,n}(0,\ldots,0)
=[\tau_0^n]^\Theta_g\) for the cusp volume.

The descendant theory underlying these volumes is governed by the
Br\'ezin--Gross--Witten (BGW) tau-function of the KdV hierarchy.
The identification with theta intersection numbers, conjectured in \cite{norbury2023new}, is established in \cite{CGG,norbury2020enumerative}. The work of
Chidambaram--Garcia-Failde--Giacchetto~\cite{CGG} places this result in
a broader theory of negative \(r\)-spin classes and derives
\(\mathcal W\)-constraints for their descendant potentials.
On the recursion side, Do and Norbury~\cite{NP1} studied the Bessel
spectral curve and proved that its partition function satisfies the
KdV hierarchy, Virasoro constraints, and a cut-and-join recursion.
Huang--Liu--Xu~\cite{HLX} extended the Stanford--Witten volume recursion
to intersections involving higher kappa classes. These results provide
the geometric and recursive framework for the shifted theta
intersection numbers used below.

Large-genus asymptotics add an analytic aspect to this framework.
Guo--Norbury--Yang--Zagier~\cite{GNYZ} proved a uniform leading
asymptotic theorem for normalized BGW numbers, studied polynomiality
of their asymptotic coefficients, and developed applications to
BGW--kappa intersections. Their uniform theorem is an essential input
here, because the auxiliary BGW endpoint in our argument has a number
of insertions growing with the genus. From the matrix-model side,
Griguolo--Papalini--Russo--Seminara
\cite[equations~(17) and~(15a)]{griguolo2024asymptotics} proposed
large-genus formulae for super Weil-Petersson volumes using resurgence
in JT supergravity. In our normalizations, the fixed-boundary
formulation of their prediction is as follows.

\begin{conjecture}[GPRS \cite{griguolo2024asymptotics}]\label{Mainconj}
(1) For every fixed \(n\geq0\) and fixed nonnegative boundary
lengths \(L_1,\ldots,L_n\), as \(g\to\infty\),
\begin{equation}\label{Vconj}
 V^\Theta_{g,n}(L_1,\ldots,L_n)
 \sim\frac{\Gamma(2g+n-\frac52)}
 {2^{n-\frac12}\pi^{\frac92-n-2g}}
 \prod_{i=1}^n\frac{\sinh(L_i/4)}{L_i/4}.
\end{equation}
(2) For every fixed \(n\geq0\) and fixed
\(\mathbf d\in\mathbb Z_{\geq0}^n\), as \(g\to\infty\),
\begin{equation}\label{Coeffconj}
 [\tau_{d_1}\cdots\tau_{d_n}]^\Theta_g
 \sim\frac{2^{-n-2|\mathbf d|}\pi^{2g+n-\frac92}}
 {\sqrt{g-|\mathbf d|-1}}(2g-3+n)!.
\end{equation}
The quotient \(\sinh(L_i/4)/(L_i/4)\) is interpreted as one when
\(L_i=0\), and products are empty when \(n=0\).
\end{conjecture}

Both parts of Conjecture~\ref{Mainconj} are proved in
Section~\ref{sec:gprs-proof}: part~(1) is
Theorem~\ref{gprs:boundary-volumes}, and part~(2) is
Theorem~\ref{gprs:coefficients}. The proof also gives uniformity of
\eqref{Vconj} when \(\sum_iL_i=o(\sqrt g)\).
Proposition~\ref{gprs:obstruction} shows why this relative asymptotic
cannot be extended to every family with \(L_i=o(g)\).

For ordinary Weil-Petersson volumes, Mirzakhani
\cite{mirzakhani2013growth} established large-genus estimates in
connection with random hyperbolic surfaces. Mirzakhani and Zograf
\cite{MZ} subsequently proved a complete expansion in inverse powers
of the genus, with polynomial dependence on the number of markings
and a universal leading constant. The following result of the first
author is the super analogue of the Mirzakhani--Zograf asymptotic
theorem: it establishes the full expansion for \(V^\Theta_{g,n}\),
while leaving the universal constant undetermined.

\begin{theorem}[{\cite[Theorem~1.5]{Huang}}]\label{superWP}
There is a universal constant \(0<C<\infty\) such that, for every
fixed integer \(s\geq1\) and \(n\geq0\),
\[
 V^\Theta_{g,n}
 =C\frac{(2g-3+n)!\pi^{2g+n}}{2^n\sqrt g}
  \left(1+\frac{d_n^1}{g}+\cdots+\frac{d_n^s}{g^s}
         +O_{n,s}(g^{-s-1})\right)
 \qquad(g\to\infty).
\]
For each \(i\geq1\), the coefficient \(d_n^i\) is a polynomial in
\(n\) of degree \(2i\), with coefficients in
\(\mathbb Q[\pi^{-1}]\), and can be computed recursively.
\end{theorem}

Our first main result identifies this constant.

\begin{theorem}\label{Mainthm}
The universal constant in Theorem~\ref{superWP} is
\[
 C=\pi^{-9/2}.
\]
Consequently, for every fixed \(n\geq0\),
\[
 V^\Theta_{g,n}
 =\frac{\Gamma(2g+n-2)\pi^{2g+n-9/2}}{2^n\sqrt g}
  \left(1+O_n(g^{-1})\right).
\]
\end{theorem}

The proof adapts the uniform-transport strategy in a previous work of the last two authors
to a family of shifted theta intersection numbers. One endpoint is
the super Weil-Petersson specialization; the other has an exact
expression in terms of BGW numbers, to which the uniform theorem of
Guo--Norbury--Yang--Zagier \cite{GNYZ} applies. Positivity, separating-convolution estimates, and
a first descendant correction with a summable remainder show that a
suitably normalized amplitude has logarithmic derivative
\(O(g^{-1})\) along the family. Transporting the endpoint amplitude
then determines \(C\). The super identities and the estimates uniform
in the interpolation parameter are established explicitly in
Sections~\ref{sec:identities}--\ref{sec:transport}.

We also compute the first four correction coefficients and establish
general bounds on their polynomial structure. Put \(q=\pi^{-1}\)
and define \(\ell_r(n)\) by the formal identity
\[
 \log\left(1+\sum_{r\geq1}d_n^r x^r\right)
 =\sum_{r\geq1}\ell_r(n)x^r.
\]

\begin{theorem}[Structure of the volume coefficients]\label{Structurethm}
For every \(r\geq1\),
\[
 \ell_r(n),d_n^r\in\mathbb Q[n,q],\qquad
 \deg_n\ell_r\leq r+1,\qquad \deg_n d_n^r=2r.
\]
The two highest coefficients in \(n\) are
\begin{equation}
\begin{aligned}
 \relax[n^{2r}]d_n^r&=\frac{(-1)^r q^{2r}}{r!},\\
 [n^{2r-1}]d_n^r&=
 \frac{(-1)^r q^{2r-2}}{(r-1)!}
 \left(\frac{2r-1}{4}-(2r+1)q^2\right).
\end{aligned}
 \label{eq:intro-structure-top-two}
\end{equation}
The following bounds on inverse-pi degree hold:
\begin{equation}
 \deg_q\ell_r\leq2r+2,\qquad
 \deg_q\bigl([n^j]\ell_r\bigr)\leq2r\quad(j\geq1),\qquad
 \deg_q d_n^r\leq2r+2\lfloor r/2\rfloor.
 \label{eq:intro-structure-q-bounds}
\end{equation}
These inverse-pi bounds are not claimed to be optimal.
\end{theorem}

The explicit formulae for \(d_n^1,\ldots,d_n^4\) appear in
\eqref{d1}--\eqref{d4}. Sections~\ref{sec:interpolation}--\ref{sec:constant}
prove Theorem~\ref{Mainthm}, and Section~\ref{sec:gprs-proof} gives
the detailed proof of Conjecture~\ref{Mainconj}.
Sections~\ref{sec:coefficient-procedure} and
\ref{sec:coefficient-structure} develop the coefficient calculation
and its structural consequences. The appendices provide an independent
KdV verification and exact and numerical checks. Unless stated
otherwise, the number of markings is fixed; constants in the
interpolation estimates may depend on this number but are uniform
in \(t\in[0,1]\).

For dependence of \(V^\Theta_{g,n}\) in $n$, we have the following theorem, which will be proved in Section \ref{sec:large-n}.
\begin{theorem}\label{polyin-n}
For fixed \(g\geq2\), the quotient
\(V^\Theta_{g,n}/n!\) is a polynomial in \(n\).
\end{theorem}

\

\noindent{\bf Disclosure of AI assistance.}
GPT-6 Astra was used as an auxiliary tool for literature search, proofreading, 
and checks of mathematical arguments.

\section{A shifted BGW interpolation}
\label{sec:interpolation}

Define, with the analytic branch equal to one at the origin,
\begin{equation}
H(s)=\prod_{j\geq1}\left(1-\frac{s}{(2j+1)^2}\right)^{-1},
\qquad f_t(s)=(1-s)H(s)^{-t}.
\label{eq:family}
\end{equation}
Euler's cosine product gives
\[
f_0(s)=1-s,\qquad
f_1(s)=\cos\left(\frac{\pi}{2}\sqrt s\right).
\]
In spectral-curve language, this interpolates between
\(y_0(z)=z^{-1}-16z\) and \(y_1(z)=\cos(2\pi z)/z\), keeping the first
positive zero at \(z=1/4\). The proof uses shifted intersection numbers
directly, so no additional spectral-curve identification is assumed.

Write \(f_t(s)=\sum_{k\geq0}c_k(t)s^k\), and put
\[
T_k(t)=-16^k(2k-1)!!c_k(t)\quad(k\geq1),
\qquad p_k=16^k(2k+1)!!.
\]
Define
\begin{equation}
\begin{aligned}
F^{(t)}_{g,n}(\dd)
&=\left\langle
 \prod_{i=1}^n\tau_{d_i}
 \exp\left(\sum_{k\geq1}T_k(t)\tau_k\right)
 \right\rangle_g^\Theta,\\
J^{(t)}_{g,n}(\dd)
&=\left(\prod_{i=1}^n p_{d_i}\right)F^{(t)}_{g,n}(\dd),
\qquad V^{(t)}_{g,n}=J^{(t)}_{g,n}(0^n).
\end{aligned}
\label{eq:shifted}
\end{equation}
The exponential has its usual exponential-generating interpretation.
Each expression is finite: a monomial in additional insertions satisfies
\[
|\dd|+\sum_{k\geq1}k m_k=g-1,
\qquad |\dd|=\sum_i d_i.
\]

At \(t=1\),
\[
T_k(1)=\frac{(-1)^{k+1}(2\pi^2)^k}{k!}.
\]
The theta forgetful-map identity relates kappa parameters \(s_j\) to
descendant shifts by
\[
\sum_{k\geq1}T_kz^k
=1-\exp\left(-\sum_{j\geq1}s_jz^j\right).
\]
Taking \(s_1=2\pi^2\) and \(s_{j>1}=0\) gives
\begin{equation}
V^{(1)}_{g,n}=V^\Theta_{g,n},\qquad
J^{(1)}_{g,n}(\dd)
=4^{|\dd|}[\tau_{d_1}\cdots\tau_{d_n}]_g^\Theta.
\label{eq:identification}
\end{equation}
For this specialization, the identity follows explicitly from
\cite[equation~(194)]{GNYZ}: multiply the formula for
\(\langle\kappa_1^m\prod_i\tau_{d_i}\rangle_g^\Theta\)
by \(a^m/m!\) and sum over \(m\). The resulting shifts are
\(T_k=(-1)^{k+1}a^k/k!\).
Thus \eqref{eq:identification} uses an exact identity, not an
asymptotic assertion about super volumes.

At \(t=0\), the only nonzero shift is \(T_1=16\).

\section{Kernel bounds and exact identities}
\label{sec:identities}

Set
\begin{equation}
\begin{aligned}
\frac{1}{f_t(s)}
&=\sum_{L\geq0}b_L(t)s^L,
& b_{-1}&=0,
& \Delta_L&=b_L-b_{L-1},\\
B_t&=(4/\pi)^t,
& M_t&=\frac{t}{4}B_t,
& \alpha_t&=\int_0^1 f_t(x^2)\,dx.
\end{aligned}
\label{eq:kernel}
\end{equation}
Since \(\Delta_L=[s^L]H(s)^t\) and \(\log H\) has positive Taylor
coefficients,
\begin{equation}
\begin{gathered}
1=b_0\leq b_L\leq B_t,\qquad \Delta_L\geq0,\\
\sum_{L\geq0}\Delta_L=B_t,\qquad
\sum_{L\geq0}L\Delta_L=M_t.
\end{gathered}
\label{eq:kernel-moments}
\end{equation}
Here \(H(1)=4/\pi\) and
\[
\frac{H'(1)}{H(1)}
=\sum_{j\geq1}\frac{1}{(2j+1)^2-1}
=\frac14.
\]
For every \(1<\rho<9\), analyticity on \(|s|\leq\rho\) implies
\(\Delta_L\leq C_\rho\rho^{-L}\), uniformly in \(t\).
All polynomial moments are consequently bounded uniformly. Also,
\begin{equation}
\alpha_0=\frac23,\qquad \alpha_1=\frac2\pi,
\qquad \frac2\pi\leq\alpha_t\leq\frac23.
\label{eq:alpha}
\end{equation}

We suppress \(t\) temporarily. Write \(D=g-1-|\dd|\), and define the
two-node expression
\[
\begin{aligned}
\mathcal R_{g,n}(u,v;\dd_{\widehat{1}})
={}&J_{g-1,n+1}(u,v,\dd_{\widehat{1}})\\
&+\sum_{\substack{h+h'=g\\
 I\sqcup I'=\{2,\ldots,n\}}}
 J_{h,|I|+1}(u,\dd_I)
 J_{h',|I'|+1}(v,\dd_{I'}).
\end{aligned}
\]
The shifted super DVV recursion, except for the initial pair
\((g,n)=(1,1)\), is
\begin{equation}
\begin{aligned}
J_{g,n}(\dd)
={}&\sum_{j=2}^n(2d_j+1)\sum_{L=0}^{D}b_L
 J_{g,n-1}(L+d_1+d_j,\dd_{\widehat{1},\widehat{j}})\\
&+8\sum_{L=0}^{D}b_L
 \sum_{\substack{u,v\geq0\\u+v=L+d_1-1}}
 \mathcal R_{g,n}(u,v;\dd_{\widehat{1}}).
\end{aligned}
\label{eq:recursion}
\end{equation}
Splittings are ordered; negative degrees, excess degrees, and
genus-zero factors mean zero. The initial value is \(J_{1,1}(0)=1/8\).

For clarity, \eqref{eq:recursion} follows directly from the
super DVV formula \cite[equation~(1.8)]{HLX}, which was firstly proved by Do and Norbury \cite{NP1}. In the shifted
exponential, merging the distinguished insertion with an additional
\(\tau_k\) contributes
\(T_k/[16^k(2k-1)!!]=-c_k\) times \(J(d_1+k,\ldots)\).
Moving these terms to the left gives
\[
\sum_{k\geq0}c_kJ(d_1+k,\ldots)=R_0(\dd),
\]
where \(R_0\) consists of the unshifted merge and node terms.
The degree-raising operator is nilpotent at fixed genus, so
multiplication by its inverse series \(1/f_t\) gives
\eqref{eq:recursion} exactly.

The theta dilaton identity and the first integrated KdV equation give
\begin{equation}
(2g-2+n)V_{g,n}
=\sum_{k\geq0}w_kJ_{g,n+1}(k,0^n),
\qquad w_k=\frac{c_k}{2k+1},
\label{eq:dilaton}
\end{equation}
and, for \(n\geq2\),
\begin{equation}
\begin{aligned}
J_{g,n}(1,0^{n-1})
={}&4V_{g-1,n+2}\\
&+24\sum_{h=1}^{g-1}\sum_{r=0}^{n-2}
 \binom{n-2}{r}V_{h,r+2}V_{g-h,n-r}.
\end{aligned}
\label{eq:kdv}
\end{equation}
In \eqref{eq:dilaton}, the extra-marking term in the dilaton identity
is \(\sum_k T_k\partial_{T_k}\); moving it to the other side uses
\(-T_k/p_k=w_k\) for \(k\geq1\), with \(w_0=1\).
In \eqref{eq:kdv}, multiplying the integrated KdV coefficients
\(1/12,1/2\) by \(p_1=48\) gives \(4,24\).

There is no undetermined integration constant in this KdV identity.
Indeed, for the logarithm \(\mathcal F\) of the BGW tau-function and
\(x=t_0\), KdV implies that
\[
\mathcal F_{t_1x}
-\frac{\hbar^2}{12}\mathcal F_{xxxx}
-\frac{\hbar^2}{2}(\mathcal F_{xx})^2
\]
is independent of \(x\). At genus \(g\), the coefficient of a monomial
with \(m\) positive-degree insertions can be nonzero only if their
degrees sum to \(g-2\). The theta dilaton identity makes each term
a multiple of \((1-x)^{-(2g+m)}\). Since \(2g+m>0\), independence
of \(x\) forces the difference to vanish. Taking additional
derivatives and then shifting the times proves \eqref{eq:kdv}.

Analyticity of \(f_t\) and \(\partial_tf_t\) also gives, for every
fixed integer \(m\geq0\),
\begin{equation}
\sup_t\sum_{k\geq0}(1+k)^m
\bigl(|w_k|+|\dot w_k|\bigr)<\infty.
\label{eq:weight-bounds}
\end{equation}

\section{The BGW endpoint and its exact amplitude}
\label{sec:endpoint}

Use \(C_{\BGW}(\dd)\) for the normalized numbers in \cite{GNYZ}, to
distinguish them from the unknown constant:
\begin{equation}
C_{\BGW}(\dd)
=\frac{2^{2g-1}\prod_{i=1}^N(2d_i+1)!!}
{\Gamma(2g-2+N)}
\langle\tau_{d_1}\cdots\tau_{d_N}\rangle_g^\Theta.
\label{eq:bgw-normalization}
\end{equation}
Definition~(6), identity~(8), and Theorem~1, equations~(19)--(20),
of \cite{GNYZ} give
\begin{equation}
C_{\BGW}(0,\dd)=C_{\BGW}(\dd),\qquad
\left|C_{\BGW}(\dd)-\pi^{-1}\right|\leq\frac Kg
\label{eq:bgw-theorem}
\end{equation}
with an absolute \(K\), without fixing the number of markings.
Set \(D_g=C_{\BGW}(1^{g-1})\). At \(t=0\),
\eqref{eq:shifted} and \eqref{eq:bgw-normalization} imply, for \(g\geq2\),
\begin{equation}
\begin{aligned}
V^{(0)}_{g,n}
&=\frac{16^{g-1}}{(g-1)!}
  \langle\tau_0^n\tau_1^{g-1}\rangle_g^\Theta\\
&=D_g\frac{16^{g-1}\Gamma(3g+n-3)}
 {2^{2g-1}3^{g-1}\Gamma(g)},
\qquad D_g=\frac1\pi+O(g^{-1}).
\end{aligned}
\label{eq:endpoint-exact}
\end{equation}

For \(n=2\), Stirling's formula yields
\[
\frac{\Gamma(3g-1)}
{\Gamma(g)\Gamma(2g-1/2)}
=\sqrt{\frac2\pi}\,3^{-3/2}
\left(\frac{27}{4}\right)^g
\bigl(1+O(g^{-1})\bigr).
\]
Consequently,
\begin{equation}
\frac{V^{(0)}_{g,2}(2/3)^{2g-1/2}}
{4^g\Gamma(2g-1/2)}
=\frac1{8\pi^{3/2}}\bigl(1+O(g^{-1})\bigr).
\label{eq:endpoint-amplitude}
\end{equation}

We also need a coarse consequence of \eqref{eq:endpoint-exact}:
for each fixed \(k\geq1\), there are \(C_k,P_k\) such that
\begin{equation}
V^{(0)}_{h,k}\leq C_k36^h(h!)^2h^{P_k},
\qquad h\geq1.
\label{eq:endpoint-bound}
\end{equation}
This follows from boundedness of \(D_h\) and Stirling's formula;
the finitely many small genera, including genus one, are absorbed
in the constant.


\section{Uniform coarse estimates}
\label{sec:coarse}

The parameter \(t\in[0,1]\) is suppressed in most formulas. A constant
may depend on the fixed marking numbers appearing in a statement, but
never on \(t\), the genus, or the descendant degrees. We use the
conventions in Section~\ref{sec:identities}: genus-zero factors and
descendants outside their dimension range are zero.

For later estimates, record the consequences of analyticity
\[
K_m:=\sup_{0\leq t\leq1}\sum_{L\geq0}(1+L)^m\Delta_L(t)<\infty,
\qquad
W_m:=\sup_{0\leq t\leq1}\sum_{k\geq0}(1+k)^m
\bigl(|w_k(t)|+|\dot w_k(t)|\bigr)<\infty
\]
for every fixed integer \(m\geq0\). Moreover, for any
\(1<\rho_0<\rho<9\), Cauchy's bound for \(H(s)^t\) gives
\begin{equation}
\sum_{L>D}(1+L)^m\Delta_L(t)
\leq C_{m,\rho_0,\rho}\rho_0^{-D}
\qquad(D\geq0).
\label{uniform:tail}
\end{equation}
Indeed, \(\Delta_L\leq C_\rho\rho^{-L}\), and
\(\sup_L(1+L)^m(\rho_0/\rho)^L<\infty\).
These bounds justify all infinite moment sums and the later removal of
finite degree cutoffs.

\begin{lemma}[Positivity, comparison, and adjacent differences]
\label{uniform:positivity}
For \(g,n\geq1\), cusp volumes are positive and
\begin{equation}
J^{(0)}_{g,n}(\dd)
\leq J^{(t)}_{g,n}(\dd)
\leq B_t^{2g-2+n}J^{(0)}_{g,n}(\dd).
\label{eq:comparison}
\end{equation}
For every nonnegative degree vector,
\begin{equation}
0\leq J_{g,n}(\dd+\mathbf e_i)
\leq J_{g,n}(\dd)\leq V_{g,n}.
\label{eq:monotonicity}
\end{equation}
Except for the initial pair \((g,n)=(1,1)\), the exact adjacent
difference in the first coordinate is the right-hand side of
\eqref{eq:recursion} with \(b_L\) replaced by \(\Delta_L\), with
the original degrees and original upper limit \(D=g-1-|\dd|\).
\end{lemma}

\begin{proof}
For the polynomial comparison, temporarily regard the initial value as
\(b_0/8\). Induction on \(2g-2+n\) in \eqref{eq:recursion} shows
that each admissible \(J_{g,n}(\dd)\) is a polynomial in the
\(b_L\) with nonnegative coefficients, homogeneous of ordinary degree
\(2g-2+n\). A merge child and a nonseparating child both have degree
\(2g-3+n\). For a separating term, the sum of the degrees of its two
factors is also \(2g-3+n\). Multiplication by \(b_L\) supplies the
remaining degree. At the actual specialization \(b_0=1\), and
\(1\leq b_L(t)\leq B_t\); evaluating the homogeneous polynomial
proves \eqref{eq:comparison}. Positivity of cusp volumes also follows
directly from the recursion. In genus one,
\(V_{1,n}=(n-1)!/8>0\). For \(g\geq2\), the nonseparating term
with \(L=1\) and \(u=v=0\) contributes
\(8b_1V_{g-1,n+1}>0\). Induction on the genus proves strict
positivity throughout the interpolation.

Let \(R_L(\dd)\) denote the merge and node expression multiplying
\(b_L\) in \eqref{eq:recursion}. Raising \(d_1\) by one changes
each child index exactly as raising \(L\) by one. The merge
multiplicities \(2d_j+1\), \(j\geq2\), do not change. Therefore
\[
\begin{aligned}
J(\dd)&=\sum_{L=0}^{D}b_LR_L(\dd),\\
J(\dd+\mathbf e_1)
&=\sum_{L=0}^{D-1}b_LR_{L+1}(\dd)
=\sum_{L=1}^{D}b_{L-1}R_L(\dd).
\end{aligned}
\]
Subtraction proves the asserted difference formula, including \(D=0\),
when the second sum is empty. All its terms are nonnegative because
\(\Delta_L\geq0\). The initial pair has values
\(J_{1,1}(0)=1/8\), \(J_{1,1}(d)=0\) for \(d>0\), so monotonicity
holds there as well. The shifted intersection definition is symmetric
in the markings. Hence the same conclusion holds in each coordinate;
iteration from \(\mathbf0\) gives \eqref{eq:monotonicity}.
\end{proof}

\begin{lemma}[One-sided bounds]
\label{uniform:one-sided}
There is a constant \(C_*\), independent of \(g,n,t\), such that
\begin{equation}
\begin{aligned}
(2g-2+n)V_{g,n}&\leq C_*V_{g,n+1}
&& (g,n\geq1),\\
4V_{g-1,n+2}&\leq V_{g,n}
&& (g\geq1,\ n\geq2).
\end{aligned}
\label{eq:one-sided}
\end{equation}
For every fixed \(k\geq1\), there are constants \(C_k,P_k\) such that
\begin{equation}
V^{(t)}_{h,k}\leq C_k B_t^{2h}36^h(h!)^2h^{P_k},\qquad h\geq1.
\label{eq:factorial-bound}
\end{equation}
\end{lemma}

\begin{proof}
Take absolute values in the exact identity \eqref{eq:dilaton} and use
\eqref{eq:monotonicity}. The sum of absolute weights is bounded by
\(W_0\), so \(C_*=W_0\) suffices. Every term on the right of
\eqref{eq:kdv} is nonnegative, and its left side is at most \(V_{g,n}\);
this gives the second inequality. Finally, combine
\eqref{eq:endpoint-bound} with \eqref{eq:comparison}. The remaining
factor \(B_t^{k-2}\) is uniformly bounded for fixed \(k\), because
\(1\leq B_t\leq4/\pi\).
\end{proof}

\begin{lemma}[Coarse separating estimate]
\label{uniform:coarse-split}
For fixed \(n\geq2\) and integers \(k,b\geq1\) with
\(k+b=n+1\),
\begin{equation}
\sum_{h=1}^{g-1}V_{h,k}V_{g-h,b}=O_n(V_{g,n}/g).
\label{eq:coarse-split}
\end{equation}
\end{lemma}

\begin{proof}
If \(k,b\geq2\), then necessarily \(n\geq3\). Set \(r=k-2\)
in \eqref{eq:kdv} with \(n-1\) markings. The relevant summand has
the positive coefficient \(24\binom{n-3}{k-2}\). Thus the sum in
\eqref{eq:coarse-split} is bounded by a constant times \(V_{g,n-1}\).
The first inequality of \eqref{eq:one-sided}, applied to \(n-1\)
markings, gives \(V_{g,n-1}=O_n(V_{g,n}/g)\).

The remaining case is \(k=1,b=n\), after interchanging the two
components if necessary. Fix \(0<\epsilon\leq1/2\), to be decreased
below. For \(h\geq\epsilon g\),
\[
V_{h,1}\leq\frac{C_*}{2h-1}V_{h,2}
\leq\frac{C_\epsilon}{g}V_{h,2}.
\]
The sum of \(V_{h,2}V_{g-h,n}\) occurs with coefficient \(24\) in
\eqref{eq:kdv} for \(n\) markings. Hence this part is
\(O_n(V_{g,n}/g)\).

For the remaining range \(1\leq h<\epsilon g\), first observe that
for \(j\geq2\),
\[
\begin{aligned}
\frac{V_{j-1,n}}{V_{j,n}}
&\leq
\frac{C_*^2}{(2j+n-4)(2j+n-3)}
\frac{V_{j-1,n+2}}{V_{j,n}}\\
&\leq\frac{C_n}{j^2}.
\end{aligned}
\]
Here both denominators are positive, including \(j=n=2\). Iterating
over \(j=g-h+1,\ldots,g\), all of which exceed \(g/2\), gives
\[
\frac{V_{g-h,n}}{V_{g,n}}\leq(4C_n/g^2)^h.
\]
By \eqref{eq:factorial-bound} and \(h!\leq h^h\), fixed polynomial
factors can be absorbed into an exponential, giving
\(V_{h,1}\leq C_0(C_1h)^{2h}\), uniformly in \(t\). The normalized
small-genus sum is therefore at most
\[
C_0\sum_{1\leq h<\epsilon g}(K h/g)^{2h}.
\]
Choose \(\epsilon\) so that \(K\epsilon\leq1/2\). For each such
\(h\),
\[
(K h/g)^{2h}
\leq K^2h^2g^{-2}4^{-(h-1)}.
\]
The geometric moment \(\sum_{h\geq1}h^2 4^{-(h-1)}\) is finite, so
this part is actually \(O(g^{-2})\). This proves the claim. Enlarging
the constant handles the finitely many small genera.
\end{proof}

\begin{proposition}[A deficit bound for every degree]
\label{uniform:deficit}
For every fixed \(n\geq2\), all \(g\geq1\), and all nonnegative
degree vectors \(\dd\),
\begin{equation}
0\leq1-\frac{J_{g,n}(\dd)}{V_{g,n}}
\leq\min\{1,C_n|\dd|^2/g\}.
\label{eq:deficit}
\end{equation}
\end{proposition}

\begin{proof}
First assume \(|\dd|\leq g-1\), so the exact adjacent-difference
formula of Lemma~\ref{uniform:positivity} applies. The two non-product
cusp ratios satisfy
\[
\frac{V_{g,n-1}}{V_{g,n}}\leq\frac{C_n}{g},\qquad
\frac{V_{g-1,n+1}}{V_{g,n}}\leq\frac{C_n}{g}.
\]
The first follows from \eqref{eq:one-sided}, even when the merge child
has one marking. For the second, when \(g\geq2\), first use the
marking inequality at genus \(g-1\) to compare
\(V_{g-1,n+1}\) to \(V_{g-1,n+2}\), then the genus inequality in
\eqref{eq:one-sided}. At \(g=1\) this child has genus zero and is zero.

By monotonicity, replace each child descendant in the difference formula
by its cusp volume. The sum of merge multiplicities is
\(\sum_{j=2}^n(2d_j+1)\leq2|\dd|+n-1\). For a node term,
the number of nonnegative pairs with \(u+v=L+d_1-1\) is exactly
\(L+d_1\); it is zero if \(L=d_1=0\). For each of the finitely many
partitions of the remaining markings, Lemma~\ref{uniform:coarse-split}
controls its separating cusp sum. Summing against \(\Delta_L\), and
using \(K_0,K_1<\infty\), therefore proves
\begin{equation}
0\leq
\frac{J_{g,n}(\dd)-J_{g,n}(\dd+\mathbf e_i)}{V_{g,n}}
\leq\frac{C_n(1+|\dd|)}{g}.
\label{uniform:adjacent-coarse}
\end{equation}
This holds in any coordinate by symmetry.

For \(m=|\dd|\geq1\), choose any coordinate-increasing path of
length \(m\) from \(\mathbf0\) to \(\dd\). Summing
\eqref{uniform:adjacent-coarse} along that path gives
\[
1-\frac{J_{g,n}(\dd)}{V_{g,n}}
\leq\frac{C_n}{g}\sum_{j=0}^{m-1}(1+j)
\leq\frac{C'_n m^2}{g}.
\]
The case \(m=0\) is exact. Monotonicity supplies the upper bound by
one. Finally, beyond the dimension, \(|\dd|\geq g\), and
\(J_{g,n}(\dd)=0\); the same inequality is automatic after taking
\(C_n\geq1\). Small genera can be covered with the same constant
because the ratios are in \([0,1]\).
\end{proof}

\begin{lemma}[Uniform marking ratio]
\label{uniform:marking}
For every fixed \(n\geq1\),
\begin{equation}
\frac{V_{g,n+1}}{V_{g,n}}
=\frac{2g-2+n}{\alpha_t}+O_n(1),
\label{eq:marking}
\end{equation}
uniformly in \(t\). In particular, for fixed \(b,n\geq1\),
\begin{equation}
\frac{V_{G,b}}{V_{G,n}}=O_{b,n}(G^{b-n})
\qquad(G\longrightarrow\infty),
\label{uniform:marking-comparison}
\end{equation}
uniformly in \(t\).
\end{lemma}

\begin{proof}
Apply Proposition~\ref{uniform:deficit} with \(n+1\geq2\) markings
to the descendants in \eqref{eq:dilaton}. Since
\(\sum_{k\geq0}w_k=\int_0^1f_t(x^2)\,dx=\alpha_t\),
\[
\left|(2g-2+n)\frac{V_{g,n}}{V_{g,n+1}}-\alpha_t\right|
\leq\frac{C_n}{g}\sum_{k\geq0}k^2|w_k|
\leq\frac{C'_n}{g}.
\]
Absolute convergence follows from \(W_2<\infty\). Since
\(\alpha_t\geq2/\pi\), the right-hand denominator is bounded away
from zero for all sufficiently large \(g\), with a threshold independent
of \(t\). Uniform inversion gives \eqref{eq:marking}.
For every fixed marking number \(m\), this formula bounds
\(V_{G,m+1}/V_{G,m}\) above and below by positive constants times
\(G\). A finite product of these ratios, or their reciprocals, proves
\eqref{uniform:marking-comparison}.
\end{proof}

\section{Refined splitting and the first descendant correction}
\label{sec:correction}

\begin{lemma}[Factorial comparison across genera]
\label{uniform:genus-comparison}
For fixed \(N\geq2\) and \(1\leq h\leq g/2\),
\begin{equation}
\frac{V_{g-h,N}}{V_{g,N}}
\leq C_N\left(\frac{\alpha_t^2}{16}\right)^h
\left(\frac{(g-h)!}{g!}\right)^2.
\label{eq:iteration}
\end{equation}
\end{lemma}

\begin{proof}
Two uses of \eqref{eq:marking} at genus \(j-1\), followed by the
one-sided genus inequality in \eqref{eq:one-sided}, give
\[
\begin{aligned}
\frac{V_{j-1,N}}{V_{j,N}}
&=\frac{V_{j-1,N}}{V_{j-1,N+1}}
  \frac{V_{j-1,N+1}}{V_{j-1,N+2}}
  \frac{V_{j-1,N+2}}{V_{j,N}}\\
&\leq
\frac{\alpha_t^2}
{4(2j+N-4)(2j+N-3)}(1+O_N(j^{-1}))\\
&\leq\frac{\alpha_t^2}{16j^2}(1+C_N/j).
\end{aligned}
\]
The last bound holds for all sufficiently large \(j\), uniformly in
\(t\). For \(h\leq g/2\), all genera in the product
\(j=g-h+1,\ldots,g\) tend to infinity with \(g\). Moreover,
\[
\prod_{j=g-h+1}^g(1+C_N/j)
\leq\exp\left(C_N\sum_{j=g-h+1}^g\frac1j\right)
\leq C'_N.
\]
Multiplying the factors \(j^{-2}\) proves the claimed factorial ratio.
The finite collection of smaller \(g\) is absorbed by a constant;
all its denominators are positive and continuous in \(t\) on the
compact interval \([0,1]\).
\end{proof}

\begin{lemma}[Factorial convolution]
\label{uniform:factorial-convolution}
For fixed integers \(H\geq1\), fixed \(P\geq0\), and \(0\leq r<2\),
\[
\sum_{h=H}^{\lfloor g/2\rfloor}
\frac{h^P r^{2h}}{\binom gh^2}=O_{H,P,r}(g^{-2H}).
\]
In particular, \(H=1\) gives an \(O(g^{-2})\) bound.
\end{lemma}

\begin{proof}
Choose a fixed integer \(K>2H+P+3\), and split the sum into
\(H\leq h\leq K\), \(K<h\leq\sqrt g\), and
\(\sqrt g<h\leq g/2\). For fixed \(h\),
\(\binom gh^{-2}=O_h(g^{-2h})\), so the first sum is
\(O(g^{-2H})\). In the middle range,
\[
\frac{h^Pr^{2h}}{\binom gh^2}
\leq h^P(rh/g)^{2h}
\leq h^P(r^2/g)^h.
\]
For \(q=r^2/g\leq1/2\), the tail
\(\sum_{h>K}h^Pq^h\) is bounded by \(C_{K,P}q^{K+1}\), hence is
\(O(g^{-K-1})=O(g^{-2H})\). In the last range,
\(\binom gh\geq(g/h)^h\geq2^h\). The remaining sum is at most
\(\sum_{h>\sqrt g}h^P(r/2)^{2h}\), which is exponentially small in
\(\sqrt g\). This proves the estimate.
\end{proof}

\begin{proposition}[Refined separating estimates]
\label{uniform:refined-split}
For fixed \(n\geq2\),
\[
\sum_{h=1}^{g-1}V_{h,k}V_{g-h,b}=O_n(V_{g,n}/g^2)
\]
in each of the following cases:
\[
k+b=n+1,\quad k,b\geq1;
\qquad\text{or}\qquad
k+b=n+2,\quad k,b\geq2.
\]
\end{proposition}

\begin{proof}
Consider first \(1\leq h\leq g/2\). By
\eqref{uniform:marking-comparison},
\(V_{g-h,b}/V_{g-h,n}=O(g^{b-n})\), because \(g-h\) is comparable
with \(g\). Apply \eqref{eq:factorial-bound} to the smaller-genus
factor and \eqref{eq:iteration} to the larger-genus factor. This gives
\begin{equation}
\begin{gathered}
\frac{V_{h,k}V_{g-h,b}}{V_{g,n}}
\leq Ch^P g^{b-n}
\frac{r_*^{2h}}{\binom gh^2},\\
r_*:=\max_t\frac{B_t\alpha_t}{\alpha_0}
=\frac{12}{\pi^2}<2.
\end{gathered}
\label{eq:convolution}
\end{equation}
To verify the parameter, note that \(36/16=\alpha_0^{-2}\), and
\[
B_t\alpha_t
=\int_0^1(1-x^2)
\left(\frac{H(1)}{H(x^2)}\right)^t\,dx
\]
is increasing in \(t\), since \(1\leq H(x^2)\leq H(1)\).
At \(t=1\), it equals \(8/\pi^2\), and division by
\(\alpha_0=2/3\) gives \(r_*=12/\pi^2\).

Lemma~\ref{uniform:factorial-convolution} bounds this half of the sum by
\(O(g^{b-n-2}V_{g,n})\). Interchanging the two components handles the
other half and gives \(O(g^{k-n-2}V_{g,n})\). In the first stated
case, \(b-n=1-k\leq0\) and \(k-n=1-b\leq0\); in the second,
\(b-n=2-k\leq0\) and \(k-n=2-b\leq0\). Both halves are therefore
\(O(V_{g,n}/g^2)\). A possible central term is harmless if counted
twice. There are no genus-zero contributions.
\end{proof}

\begin{lemma}[Leading genus ratio and child ratios]
\label{uniform:genus-ratio}
For fixed \(n\geq2\),
\begin{equation}
\frac{4V_{g-1,n+2}}{V_{g,n}}=1+O_n(g^{-1}).
\label{eq:genus}
\end{equation}
For fixed \(n\geq3\),
\begin{equation}
\frac{V_{g,n-1}}{V_{g,n}}
=\frac{\alpha_t}{2g}+O_n(g^{-2}),\qquad
\frac{V_{g-1,n+1}}{V_{g,n}}
=\frac{\alpha_t}{8g}+O_n(g^{-2}).
\label{uniform:child-ratios}
\end{equation}
\end{lemma}

\begin{proof}
Divide \eqref{eq:kdv} by \(V_{g,n}\). Its left side is
\(1+O_n(g^{-1})\) by Proposition~\ref{uniform:deficit}, and the
separating contribution is \(O_n(g^{-2})\) by
Proposition~\ref{uniform:refined-split}. This proves \eqref{eq:genus}.
The first formula in \eqref{uniform:child-ratios} follows by uniformly
inverting \eqref{eq:marking} with \(n-1\) markings. For the second,
apply \eqref{eq:genus} with \(n-1\geq2\) and multiply by the first:
\[
\frac{V_{g-1,n+1}}{V_{g,n}}
=\frac{V_{g-1,n+1}}{V_{g,n-1}}
 \frac{V_{g,n-1}}{V_{g,n}}
=\frac14(1+O_n(g^{-1}))
 \left(\frac{\alpha_t}{2g}+O_n(g^{-2})\right).
\]
Only the one-sided genus inequality was used in
Lemma~\ref{uniform:genus-comparison}; the two-sided asymptotic is
obtained here, after the refined separating estimates. Thus the
dependency is noncircular.
\end{proof}

\begin{proposition}[First descendant correction with a summable remainder]
\label{uniform:first-correction}
For fixed \(n\geq3\), uniformly for \(t\in[0,1]\) and integers
\(k\geq0\),
\begin{equation}
\begin{aligned}
\frac{J_{g,n}(k,0^{n-1})}{V_{g,n}}
&=1-\frac{E_{n,k}(t)}g
+O_n\left(\frac{(k+1)^4}{g^2}\right),\\
E_{n,k}(t)
&=\frac{\alpha_tB_t}{2}
\left[k^2+(n-2+t/2)k\right].
\end{aligned}
\label{eq:first-correction}
\end{equation}
\end{proposition}

\begin{proof}
First consider an adjacent difference with \(0\leq d\leq\sqrt g\),
and write \(D=g-1-d\). Denote its normalized value by
\[
\mathcal D_{g,n}(d):=
\frac{J_{g,n}(d,0^{n-1})-J_{g,n}(d+1,0^{n-1})}{V_{g,n}}.
\]
The exact difference formula is
\[
\begin{aligned}
\mathcal D_{g,n}(d)
={}&(n-1)\sum_{L=0}^{D}\Delta_L
 \frac{J_{g,n-1}(L+d,0^{n-2})}{V_{g,n}}\\
&+8\sum_{L=0}^{D}\Delta_L
 \sum_{u+v=L+d-1}
 \frac{J_{g-1,n+1}(u,v,0^{n-1})}{V_{g,n}}
+\mathcal S_{g,n}(d),
\end{aligned}
\]
where \(\mathcal S_{g,n}(d)\geq0\) is the separating contribution.
The latter is bounded by cusp products, whose total marking number is
\(n+1\). Proposition~\ref{uniform:refined-split} therefore gives
\[
0\leq\mathcal S_{g,n}(d)
\leq\frac{C_n}{g^2}\sum_{L\geq0}\Delta_L(L+d)
\leq\frac{C'_n(d+1)}{g^2}.
\]
The finite number of partitions of the \(n-1\) remaining markings is
absorbed into \(C_n\).

In the merge term, the child has \(n-1\geq2\) markings and degree
\(L+d\), so Proposition~\ref{uniform:deficit} and
\eqref{uniform:child-ratios} imply
\[
\frac{J_{g,n-1}(L+d,0^{n-2})}{V_{g,n}}
=\frac{\alpha_t}{2g}
+O_n\left(\frac{(L+d+1)^2}{g^2}\right).
\]
For a nonseparating child, the total degree is \(L+d-1\), and the
genus is \(g-1\). The same deficit bound, with \(g-1\) in place of
\(g\), gives for sufficiently large \(g\)
\[
\frac{J_{g-1,n+1}(u,v,0^{n-1})}{V_{g,n}}
=\frac{\alpha_t}{8g}
+O_n\left(\frac{(L+d+1)^2}{g^2}\right).
\]
There are exactly \(L+d\) such pairs. After summation, the merge
error is bounded by \(C_n(d+1)^2/g^2\), and the nonseparating error
by \(C_n(d+1)^3/g^2\), because \(K_2,K_3<\infty\).

For \(g\) sufficiently large, \(d\leq\sqrt g\) ensures
\(D\geq g/2\). Equation~\eqref{uniform:tail} permits the two main
\(L\)-sums to be extended to infinity with an exponentially small
error, uniformly in \(t\) and \(d\leq\sqrt g\). Using
\(\sum\Delta_L=B_t\), \(\sum L\Delta_L=M_t\), the adjacent
difference consequently equals
\begin{equation}
\frac{\alpha_t}{g}
\left(\frac{n-1}{2}B_t+M_t+dB_t\right)
+O_n\left(\frac{(d+1)^3}{g^2}\right).
\label{uniform:adjacent-first}
\end{equation}
There is no exceptional \(d=0\) correction: the merge degree is
\(L+d\), and the node-pair count is \(L+d\), whose value at
\(L=d=0\) is already zero.

For \(1\leq k\leq\sqrt g\), sum \eqref{uniform:adjacent-first}
over \(d=0,\ldots,k-1\). The main term is
\[
\frac{\alpha_tB_t}{2g}
\left[k^2+\left(n-2+\frac{2M_t}{B_t}\right)k\right],
\]
and \(\sum_{d=0}^{k-1}(d+1)^3=O((k+1)^4)\). Since
\(2M_t/B_t=t/2\), this proves \eqref{eq:first-correction} in that
range; at \(k=0\), the equality is exact.

For \(k>\sqrt g\), monotonicity places the descendant ratio in
\([0,1]\), and the explicit expression gives
\(|E_{n,k}(t)|\leq C_n(k+1)^2\). Hence
\[
\left|\frac{J_{g,n}(k,0^{n-1})}{V_{g,n}}
-1+\frac{E_{n,k}(t)}g\right|
\leq1+\frac{C_n(k+1)^2}{g}
\leq C'_n\frac{(k+1)^4}{g^2}.
\]
This argument also applies when the degree is beyond the dimension.
Finitely many small genera are handled by increasing the constant.
\end{proof}

At \(t=1\), \(\alpha_1B_1/2=4/\pi^2\), so the correction becomes
\[
E_{n,k}(1)=\frac{4k(k+n-3/2)}{\pi^2},
\]
in agreement withthe first coefficients in \cite[equation~(3.12)]{Huang}. The all-degree remainder, rather than just
the formula for a fixed \(k\), is needed to sum against infinitely many
dilaton weights in the next section.

\section{Transport of the amplitude}
\label{sec:transport}

\begin{lemma}[Dilaton-weight moments]
\label{uniform:weight-moments}
Writing \(\alpha=\alpha_t\), \(B=B_t\), and
\(\chi_k=-\dot w_k\), one has
\begin{equation}
S_0:=\sum_{k\geq0}w_k=\alpha,\qquad
S_1:=\sum_{k\geq0}kw_k=-\frac\alpha2,\qquad
S_2:=\sum_{k\geq0}k^2w_k=\frac\alpha4-\frac1{2B}.
\label{eq:weight-moments}
\end{equation}
Moreover,
\[
\sum_{k\geq0}\chi_k=-\dot\alpha,\qquad
\sum_{k\geq0}k\chi_k=\frac{\dot\alpha}{2},\qquad
\sum_{k\geq0}k^2\chi_k
=-\frac{\dot\alpha}{4}-\frac{\dot B}{2B^2}.
\]
\end{lemma}

\begin{proof}
Termwise integration of \(f_t(x^2)=\sum c_kx^{2k}\) gives
\(S_0=\alpha\). Since \(w_k=c_k/(2k+1)\),
\[
2S_1=\sum_{k\geq0}c_k-\sum_{k\geq0}\frac{c_k}{2k+1}
=f_t(1)-\alpha=-\alpha.
\]
Similarly,
\[
2S_2=\sum_{k\geq0}kc_k-S_1=f_t'(1)-S_1.
\]
The factorization \(f_t(s)=(1-s)H(s)^{-t}\) gives
\(f_t'(1)=-H(1)^{-t}=-1/B\), proving the last identity.
All series and their \(t\)-derivatives converge absolutely and
uniformly by \(W_2<\infty\). Differentiating these identities and
using \(\chi_k=-\dot w_k\) proves the remaining formulas.
\end{proof}

\begin{proposition}[Uniform transport]
\label{uniform:transport}
For every fixed \(n\geq2\), define
\begin{equation}
\mathcal B_{g,n}(t)=
\frac{V^{(t)}_{g,n}\alpha_t^{2g+n-5/2}}
{4^g\Gamma(2g+n-5/2)\sqrt{B_t}}.
\label{eq:amplitude}
\end{equation}
Then
\[
\sup_{0\leq t\leq1}
\left|\partial_t\log\mathcal B_{g,n}(t)\right|
=O_n(g^{-1}),
\]
and consequently
\begin{equation}
\frac{\mathcal B_{g,n}(1)}{\mathcal B_{g,n}(0)}
=1+O_n(g^{-1}).
\label{eq:transport-integrated}
\end{equation}
\end{proposition}

\begin{proof}
Put
\[
N=2g-2+n,\qquad
r_k=\frac{J_{g,n+1}(k,0^n)}{V_{g,n+1}},\qquad
\chi_k=-\dot w_k.
\]
Differentiating the finite shifted-intersection polynomial
\eqref{eq:shifted} introduces one additional insertion. Since
\(\dot T_k/p_k=-\dot w_k\), and \(\dot w_0=0\),
\[
\partial_t\log V_{g,n}
=\frac{V_{g,n+1}}{V_{g,n}}\sum_{k\geq0}\chi_kr_k.
\]
The exact dilaton identity gives
\(N V_{g,n}/V_{g,n+1}=\sum_{k\geq0}w_kr_k\). These two formulas
combine to give
\begin{equation}
\partial_t\log V_{g,n}
=N\frac{\sum_{k\geq0}\chi_kr_k}{\sum_{k\geq0}w_kr_k}.
\label{eq:variation}
\end{equation}

The inserted descendant has \(n+1\geq3\) markings, so
Proposition~\ref{uniform:first-correction} applies. Set
\[
\beta=n-1+t/2,\qquad E_k=\frac{\alpha B}{2}(k^2+\beta k),
\qquad T=\sum_{k\geq0}w_kE_k,\qquad
Q=\sum_{k\geq0}\chi_kE_k.
\]
The all-degree remainder in \eqref{eq:first-correction} is bounded by
\(C_n(k+1)^4/g^2\). Therefore, using \(W_4<\infty\),
\[
\sum_{k\geq0}w_kr_k=\alpha-\frac Tg+O_n(g^{-2}),\qquad
\sum_{k\geq0}\chi_kr_k=-\dot\alpha-\frac Qg+O_n(g^{-2}),
\]
with uniform errors. This is the step that would not follow from an
asymptotic statement only for each fixed insertion degree.

Lemma~\ref{uniform:weight-moments} gives the explicit contractions
\begin{equation}
\begin{aligned}
T&=\frac{\alpha B}{2}(S_2+\beta S_1)
 =\frac{\alpha^2B}{8}(1-2\beta)-\frac\alpha4,\\
Q&=\frac{\alpha B\dot\alpha}{8}(2\beta-1)
 -\frac{\alpha\dot B}{4B},\\
Q+\frac{\dot\alpha}{\alpha}T
&=-\frac{\dot\alpha}{4}-\frac{\alpha\dot B}{4B}.
\end{aligned}
\label{eq:moment-cancellation}
\end{equation}
All of \(\alpha,B,\dot\alpha,\dot B,T,Q\) are bounded uniformly in
\(t\), and \(\alpha\geq2/\pi\). Thus the denominator in
\eqref{eq:variation} is bounded below by \(1/\pi\) once \(g\) is
sufficiently large, uniformly in \(t\). Uniform division gives
\[
\frac{-\dot\alpha-Q/g+O_n(g^{-2})}
{\alpha-T/g+O_n(g^{-2})}
=-\frac{\dot\alpha}{\alpha}
-\frac1{\alpha g}
\left(Q+\frac{\dot\alpha}{\alpha}T\right)
+O_n(g^{-2}).
\]
Multiply by \(N\), use \(N/g=2+O_n(g^{-1})\), and substitute
\eqref{eq:moment-cancellation}. The result is
\begin{equation}
\partial_t\log V_{g,n}
=-\left(2g+n-\frac52\right)\frac{\dot\alpha_t}{\alpha_t}
+\frac{\dot B_t}{2B_t}+O_n(g^{-1}).
\label{eq:transport}
\end{equation}
Taking the logarithmic derivative of \eqref{eq:amplitude} cancels
the two displayed main terms and proves the uniform derivative bound.

Finally, positivity of the volumes makes the logarithm well-defined.
For a uniform constant \(C_n\), integration over \([0,1]\) gives
\[
\left|\log\frac{\mathcal B_{g,n}(1)}{\mathcal B_{g,n}(0)}\right|
\leq\frac{C_n}{g}.
\]
Exponentiation yields
\(e^{-C_n/g}\leq\mathcal B_{g,n}(1)/\mathcal B_{g,n}(0)
\leq e^{C_n/g}\), proving \eqref{eq:transport-integrated}.
\end{proof}

Only a parent with \(n\geq2\) was used in this argument. Its inserted
descendants have at least three markings, which is exactly the range
needed for the uniform first-correction estimate. The normalization at
\(n=2\), combined withindependence of its leading
constant from \(n\) (see Theorem \ref{superWP}), is therefore sufficient for the final determination
of \(C\).

\section{Extracting the constant}
\label{sec:constant}

By Theorem \ref{superWP}, it suffices to get $C$ for $n=2$.

By \eqref{eq:endpoint-amplitude},
\[
\lim_{g\to\infty}\mathcal B_{g,2}(0)
=\frac1{8\pi^{3/2}}.
\]
At the other endpoint, use \eqref{superint}-\eqref{VTheta2}, \(\alpha_1=2/\pi\),
\(B_1=4/\pi\), and
\(\Gamma(2g)/\Gamma(2g-1/2)\sim\sqrt{2g}\). This gives
\[
\lim_{g\to\infty}\mathcal B_{g,2}(1)
=\frac{C\pi^3}{8}.
\]
Equation~\eqref{eq:transport-integrated} equates these limits. Hence
\[
\frac{C\pi^3}{8}=\frac1{8\pi^{3/2}},
\qquad C=\pi^{-9/2}.
\]
The independence of \(n\) extends this identification to
every fixed \(n\geq0\) in Theorem \ref{superWP}. This proves the Theorem \ref{Mainthm}.

\section{Proof of the GPRS conjecture and the boundary-length range}
\label{sec:gprs-proof}

We now deduce both assertions of Conjecture~\ref{Mainconj} from
Theorem~\ref{Mainthm} and the uniform descendant estimates. The following proofs are essentially due to the first author \cite[Theorem 1.7]{Huang}. The number
of markings is fixed throughout this section. In the volume assertion,
the boundary lengths are nonnegative. We prove the assertion for fixed
lengths, and more generally uniformly when their sum is
\(o(\sqrt g)\). We also explain why a relative asymptotic of the stated
form cannot hold for every family of lengths satisfying \(L_i=o(g)\).

We specialize to \(t=1\), and write
\[
 V_{g,n}=V^\Theta_{g,n},\qquad
 J_{g,n}(\dd)=J^{(1)}_{g,n}(\dd),\qquad
 R_{g,n}(\dd)=\frac{J_{g,n}(\dd)}{V_{g,n}}.
\]
Descendants whose total degree exceeds \(g-1\) are assigned the value
zero. By \eqref{eq:identification},
\begin{equation}
 J_{g,n}(\dd)=4^{|\dd|}
 [\tau_{d_1}\cdots\tau_{d_n}]^\Theta_g.
 \label{eq:gprs-normalization}
\end{equation}

\subsection{The one-marking case of the descendant estimate}

Proposition~\ref{uniform:deficit} is stated for at least two markings.
The following argument supplies the one-marking case, including the
uniformity in the descendant degree that is needed for varying lengths.

\begin{lemma}[Dilaton identity with retained descendants]
\label{gprs:retained-dilaton}
For \(g,n\geq1\) and every nonnegative degree vector \(\dd\),
\begin{equation}
 (2g-2+n)J_{g,n}(\dd)
 =\sum_{k\geq0}w_kJ_{g,n+1}(\dd,k),
 \qquad
 w_k=\frac{(-1)^k(\pi/2)^{2k}}{(2k+1)!}.
 \label{eq:gprs-retained-dilaton}
\end{equation}
At each fixed genus the sum on the right is finite.
\end{lemma}

\begin{proof}
The theta dilaton identity is
\[
 \left\langle\tau_0\prod_{i=1}^N\tau_{a_i}\right\rangle_g^\Theta
 =(2g-2+N)
 \left\langle\prod_{i=1}^N\tau_{a_i}\right\rangle_g^\Theta.
\]
Indeed, the relation
\(\Theta_{g,N+1}=\psi_{N+1}\pi^*\Theta_{g,N}\) reduces this to
\(\pi_*\psi_{N+1}=2g-2+N\). The differences between the retained
psi classes and their pullbacks are supported on the divisors where a
retained marking and the new marking lie on a three-pointed rational
tail. The restriction of \(\psi_{N+1}\) to each such divisor is zero,
so these pullback corrections do not contribute.

Consider a term in the shifted exponential with \(m_k\) additional
insertions of degree \(k\), and let \(M=\sum_km_k\). Applying the
identity above before summing the shifted series multiplies this term
by \(2g-2+n+M\). Consequently, regarding the shifts temporarily as
independent variables and using the normalization in
\eqref{eq:shifted},
\[
 J_{g,n+1}(\dd,0)
 =(2g-2+n)J_{g,n}(\dd)
   +\sum_{k\geq1}T_k\frac{\partial}{\partial T_k}J_{g,n}(\dd).
\]
Differentiation inserts one descendant, hence
\[
 \frac{\partial}{\partial T_k}J_{g,n}(\dd)
 =\frac{J_{g,n+1}(\dd,k)}{p_k}.
\]
Move the last sum to the other side and use \(w_0=1\) and
\(-T_k/p_k=w_k\) for \(k\geq1\). This proves
\eqref{eq:gprs-retained-dilaton}. All these operations are on finite
polynomials: a nonzero summand satisfies
\(|\dd|+\sum_k k m_k=g-1\).
\end{proof}

\begin{lemma}[Deficit estimate for every positive marking number]
\label{gprs:all-markings-deficit}
For every fixed \(n\geq1\), there is a constant \(C_n\) such that,
for all \(g\geq1\) and all nonnegative degree vectors \(\dd\),
\begin{equation}
 0\leq1-R_{g,n}(\dd)
 \leq\min\left\{1,\frac{C_n|\dd|^2}{g}\right\}.
 \label{eq:gprs-deficit}
\end{equation}
\end{lemma}

\begin{proof}
For \(n\geq2\), this is Proposition~\ref{uniform:deficit}. To treat
\(n=1\), put
\[
 A_d=\sum_{k\geq0}w_kR_{g,2}(d,k).
\]
Equation~\eqref{eq:gprs-retained-dilaton}, first with degree \(d\)
and then with degree zero, gives
\[
 R_{g,1}(d)=\frac{A_d}{A_0},\qquad
 A_0=(2g-1)\frac{V_{g,1}}{V_{g,2}}>0.
\]
Since \(\sum_{k\geq0}w_k=2/\pi\), the two-marking deficit estimate
implies
\[
 \left|A_0-\frac2\pi\right|
 \leq\frac{C_2}{g}\sum_{k\geq0}|w_k|k^2
 =O(g^{-1}).
\]
Thus \(A_0\geq1/\pi\) for all sufficiently large \(g\).
Monotonicity in Lemma~\ref{uniform:positivity} gives
\[
 0\leq R_{g,2}(0,k)-R_{g,2}(d,k)
 \leq1-R_{g,2}(d,k)
 \leq\frac{C_2(d+k)^2}{g}.
\]
For \(d\geq1\), it follows that
\[
 |A_0-A_d|
 \leq\frac{C_2}{g}\sum_{k\geq0}|w_k|(d+k)^2
 \leq\frac{C d^2}{g}\sum_{k\geq0}|w_k|(1+k)^2
 \leq\frac{C'd^2}{g}.
\]
The absolute weight moments are finite by \eqref{eq:weight-bounds}.
Division by \(A_0\), together with
\(0\leq R_{g,1}(d)\leq1\), proves the required estimate for large
\(g\). The case \(d=0\) is exact. Increasing the constant handles
the finitely many smaller genera, since \(d^2\geq1\) when \(d>0\).
The argument also includes degrees outside the dimension range.
\end{proof}

\subsection{Fixed descendant degrees}

\begin{theorem}[The coefficient assertion of the GPRS conjecture]
\label{gprs:coefficients}
For fixed \(n\geq1\) and fixed \(\dd\in\mathbb Z_{\geq0}^n\), put
\(m=|\dd|\). Then, as \(g\to\infty\),
\begin{equation}
 [\tau_{d_1}\cdots\tau_{d_n}]^\Theta_g
 =\frac{2^{-n-2m}\pi^{2g+n-9/2}(2g+n-3)!}{\sqrt{g-m-1}}
  \left(1+O_{n,\dd}(g^{-1})\right).
 \label{eq:gprs-coefficients}
\end{equation}
The same statement holds for \(n=0\), with the empty degree vector and
\(m=0\). In particular, this proves Conjecture~\ref{Mainconj}(2).
\end{theorem}

\begin{proof}
By \eqref{eq:gprs-deficit}, for fixed \(\dd\),
\(R_{g,n}(\dd)=1+O_{n,\dd}(g^{-1})\). Combining
\eqref{eq:gprs-normalization} with Theorem~\ref{Mainthm} therefore gives
\[
 [\tau_{d_1}\cdots\tau_{d_n}]^\Theta_g
 =\frac{2^{-n-2m}\pi^{2g+n-9/2}\Gamma(2g+n-2)}{\sqrt g}
  \left(1+O_{n,\dd}(g^{-1})\right).
\]
Here \(\Gamma(2g+n-2)=(2g+n-3)!\). For \(g>m+1\),
\[
 \frac1{\sqrt g}
 =\frac1{\sqrt{g-m-1}}
   \left(1-\frac{m+1}{g}\right)^{1/2}
 =\frac1{\sqrt{g-m-1}}\left(1+O_m(g^{-1})\right).
\]
This proves \eqref{eq:gprs-coefficients}. For \(n=0\), the left-hand
side is the cusp volume, so the same calculation starts directly from
Theorem~\ref{Mainthm}; no descendant estimate is needed.
\end{proof}

\subsection{Boundary lengths and a quantitative summation estimate}

Define
\[
 S(x)=\frac{\sinh x}{x}\quad(x>0),\qquad S(0)=1,
 \qquad
 \mathcal A_{g,n}=\frac{\Gamma(2g+n-5/2)\pi^{2g+n-9/2}}{2^{n-1/2}}.
\]
The notation \(\mathcal A_{g,n}\) in this subsection denotes the leading
prefactor in Conjecture~\ref{Mainconj}(1).

\begin{theorem}[The volume assertion and its uniform range]
\label{gprs:boundary-volumes}
For every fixed \(n\geq1\), all \(g\geq1\), and all
\(L_1,\ldots,L_n\geq0\),
\begin{equation}
 0\leq1-
 \frac{V^\Theta_{g,n}(L_1,\ldots,L_n)}
      {V_{g,n}\prod_{i=1}^n S(L_i/4)}
 \leq\frac{C_n}{g}\left(\sum_{i=1}^nL_i\right)^2.
 \label{eq:gprs-boundary-deficit}
\end{equation}
Consequently, as \(g\to\infty\),
\begin{equation}
 \frac{V^\Theta_{g,n}(L_1,\ldots,L_n)}
      {\mathcal A_{g,n}\prod_{i=1}^nS(L_i/4)}
 =1+O_n\!\left(\frac{1+(\sum_iL_i)^2}{g}\right),
 \label{eq:gprs-boundary-asymptotic}
\end{equation}
with an implied constant independent of the nonnegative lengths.
In particular, Conjecture~\ref{Mainconj}(1) holds for fixed lengths,
and uniformly on families satisfying \(\sum_iL_i=o(\sqrt g)\).
For \(n=0\), its assertion is the cusp-volume asymptotic with the
empty product equal to one.
\end{theorem}

\begin{proof}
Put \(x_i=L_i/4\). Replacing \(L_i\) by \(L_i/2\) in
\eqref{VTheta2}, and then using \eqref{eq:gprs-normalization}, gives
the exact identity
\begin{equation}
 \frac{V^\Theta_{g,n}(L_1,\ldots,L_n)}{V_{g,n}}
 =\sum_{\dd\in\mathbb Z_{\geq0}^n}
   R_{g,n}(\dd)\prod_{i=1}^n\frac{x_i^{2d_i}}{(2d_i+1)!}.
 \label{eq:gprs-volume-sum}
\end{equation}
The extension from the finite degree range to all \(\dd\) is exact
because \(R_{g,n}(\dd)=0\) beyond that range. In particular, both
factors of two in the boundary rescaling are accounted for: the factor
\(4^{-|\dd|}\) from the descendants combines with
\(\prod_i(L_i/2)^{2d_i}\) to give \(\prod_ix_i^{2d_i}\).

For fixed \(\mathbf L\), the terms in \eqref{eq:gprs-volume-sum}
are dominated, independently of \(g\), by the summable series
\[
 \sum_{\dd\geq0}\prod_i\frac{x_i^{2d_i}}{(2d_i+1)!}
 =\prod_i S(x_i).
\]
Moreover, \eqref{eq:gprs-deficit} implies
\(R_{g,n}(\dd)\to1\) for every fixed \(\dd\). Dominated
convergence already proves the fixed-length limit
\[
 \frac{V^\Theta_{g,n}(\mathbf L)}{V_{g,n}}
 \longrightarrow\prod_iS(x_i).
\]
To obtain a uniform estimate, we keep the degree dependence of the
deficit instead of passing to the limit term by term.

For \(x>0\), introduce a nonnegative integer-valued random variable
\(K_x\) with distribution
\[
 \Pr(K_x=k)=\frac{x^{2k}}{(2k+1)!S(x)},\qquad k\geq0.
\]
For \(x=0\), set \(K_0=0\). The probabilities sum to one by the
Taylor series for \(S\). Let \(K_1,\ldots,K_n\) be independent
variables with parameters \(x_1,\ldots,x_n\), respectively. Dividing
\eqref{eq:gprs-volume-sum} by \(\prod_iS(x_i)\) and using
\eqref{eq:gprs-deficit} gives
\begin{equation}
 0\leq1-
 \frac{V^\Theta_{g,n}(\mathbf L)}{V_{g,n}\prod_iS(x_i)}
 \leq\frac{C_n}{g}\,\mathbb E\left(\sum_iK_i\right)^2.
 \label{eq:gprs-expectation-bound}
\end{equation}
All the moments below are finite, so summing the deficit is justified
absolutely, including the degrees beyond the dimension range.

Let \(\mathcal D=(x/2)\partial_x\). Since
\(\mathcal D(x^{2k})=kx^{2k}\), direct differentiation gives
\begin{equation}
 \mathbb E K_x=\frac{\mathcal D S(x)}{S(x)}
   =\frac{x\coth x-1}{2},\qquad
 \mathbb E K_x^2=\frac{\mathcal D^2S(x)}{S(x)}
   =\frac{x^2-x\coth x+1}{4}.
 \label{eq:gprs-moments}
\end{equation}
These formulas extend continuously to \(x=0\). For \(x\geq0\),
\(1\leq x\coth x\leq1+x\), with the value at zero interpreted
by continuity. The lower bound follows from \(\tanh x\leq x\);
the upper bound follows from
\(x\coth x=x+2x/(e^{2x}-1)\leq x+1\). Thus
\[
 \mathbb E K_x\leq\frac{x}{2},\qquad
 \mathbb E K_x^2\leq\frac{x^2}{4}.
\]
By independence,
\[
 \mathbb E\left(\sum_iK_i\right)^2
 =\sum_i\mathbb E K_i^2
   +2\sum_{i<j}\mathbb E K_i\,\mathbb E K_j
 \leq\frac14\left(\sum_ix_i\right)^2
 =\frac1{64}\left(\sum_iL_i\right)^2.
\]
Substitution into \eqref{eq:gprs-expectation-bound}, with a change of
constant, proves \eqref{eq:gprs-boundary-deficit}.

It remains to identify the prefactor. With
\(z=2g+n-5/2\), the gamma-ratio expansion gives
\[
 \frac{\Gamma(z+1/2)}{\Gamma(z)}
 =z^{1/2}\left(1+O(z^{-1})\right)
 =\sqrt{2g}\left(1+O_n(g^{-1})\right).
\]
Consequently Theorem~\ref{Mainthm} becomes
\begin{equation}
 \frac{V_{g,n}}{\mathcal A_{g,n}}=1+O_n(g^{-1}).
 \label{eq:gprs-cusp-prefactor}
\end{equation}
For sufficiently large \(g\), this ratio is bounded by a constant
depending only on \(n\). Multiplying
\eqref{eq:gprs-boundary-deficit} by this ratio proves
\eqref{eq:gprs-boundary-asymptotic}. If \(\sum_iL_i=o(\sqrt g)\),
its error tends to zero. The same conclusion is uniform on every set
\(\sum_iL_i\leq\epsilon_g\sqrt g\) with \(\epsilon_g\to0\).
For \(n=0\), equation~\eqref{eq:gprs-cusp-prefactor} alone proves
the assertion, for the stable genera \(g\geq2\).
\end{proof}

\subsection{Why the range cannot include every \texorpdfstring{\(L_i=o(g)\)}{Li=o(g)}}

The next proposition distinguishes a fixed-length assertion from an
assertion uniform over all sublinear boundary lengths. It also shows
that the scale \(\sqrt g\) in the sufficient condition above reflects
a nonvanishing correction, rather than only a loss in the estimate.

\begin{proposition}[An obstruction at the square-root scale]
\label{gprs:obstruction}
For every fixed \(n\geq3\), there exists \(c_n>0\) such that, for
every fixed \(0<c<c_n\),
\begin{equation}
 \limsup_{g\to\infty}
 \frac{V^\Theta_{g,n}(c\sqrt g,0,\ldots,0)}
      {V_{g,n}S(c\sqrt g/4)}<1.
 \label{eq:gprs-square-root-obstruction}
\end{equation}
In particular, the relative asymptotic in
Conjecture~\ref{Mainconj}(1) does not hold for every choice of
\(L_i=o(g)\).
\end{proposition}

\begin{proof}
Let \(x=c\sqrt g/4\), let \(K=K_x\), and set
\[
 Q_g(c)=
 \frac{V^\Theta_{g,n}(c\sqrt g,0,\ldots,0)}{V_{g,n}S(x)}.
\]
Equation~\eqref{eq:gprs-volume-sum} says that
\(Q_g(c)=\mathbb E R_{g,n}(K,0^{n-1})\).
At \(t=1\), Proposition~\ref{uniform:first-correction} gives
\[
 R_{g,n}(k,0^{n-1})
 =1-\frac{4k(k+n-3/2)}{\pi^2g}
   +O_n\left(\frac{(k+1)^4}{g^2}\right),
\]
uniformly over all \(k\geq0\), including degrees beyond the
dimension range. Averaging this estimate is therefore legitimate and
gives
\begin{equation}
 1-Q_g(c)
 =\frac4{\pi^2g}\mathbb E[K(K+n-3/2)]
   +O_n\left(\frac{\mathbb E(K+1)^4}{g^2}\right).
 \label{eq:gprs-averaged-correction}
\end{equation}
The first two moments in \eqref{eq:gprs-moments} imply
\[
 \mathbb E[K(K+n-3/2)]
 =\frac{x^2+(2n-4)(x\coth x-1)}4.
\]
For completeness, the fourth moment needed to control the error is
\begin{equation}
 \mathbb E(K+1)^4
 =\frac{x^4+10x^3\coth x+25x^2+15x\coth x+1}{16}.
 \label{eq:gprs-fourth-moment}
\end{equation}
One way to verify this is to put \(N=2K+1\). Its distribution is
\(\Pr(N=2k+1)=x^{2k+1}/((2k+1)!\sinh x)\). Differentiating
\(\sinh x\) shows that its falling factorial moments are \(x^j\)
for even \(j\), and \(x^j\coth x\) for odd \(j\). Expanding
\((N+1)^4\) in falling factorials and using \(K+1=(N+1)/2\)
gives \eqref{eq:gprs-fourth-moment}.

For each fixed \(c>0\), we have \(x\to\infty\), and hence
\[
 \frac4{\pi^2g}\mathbb E[K(K+n-3/2)]
 \longrightarrow\frac{c^2}{16\pi^2},\qquad
 \frac{\mathbb E(K+1)^4}{g^2}
 \longrightarrow\frac{c^4}{4096}.
\]
The constant implicit in \eqref{eq:gprs-averaged-correction} depends
only on \(n\), not on \(c\). In particular, there is \(M_n<\infty\)
such that, for every fixed \(0<c\leq1\),
\[
 \limsup_{g\to\infty}
 \left|1-Q_g(c)-\frac{c^2}{16\pi^2}\right|
 \leq M_nc^4.
\]
Equivalently, the averaged correction is
\(c^2/(16\pi^2)+O_n(c^4)+o_{n,c}(1)\), with the
\(O_n(c^4)\) constant independent of \(c\in(0,1]\).
Choose \(c_n\leq1\) so small that
\(M_nc_n^2\leq1/(32\pi^2)\). Then, for \(0<c<c_n\),
\[
 \liminf_{g\to\infty}(1-Q_g(c))
 \geq\frac{c^2}{32\pi^2}>0,
\]
which proves \eqref{eq:gprs-square-root-obstruction}.
Finally, \eqref{eq:gprs-cusp-prefactor} shows that replacing the cusp
volume in its denominator by \(\mathcal A_{g,n}\) does not change the limsup.
The conjectured relative asymptotic would require this ratio to tend
to one, whereas \(c\sqrt g=o(g)\). This proves the last assertion.
\end{proof}

\section{Computing the first four correction coefficients}
\label{sec:coefficient-procedure}

We specialize throughout this section to $t=1$ and suppress the
superscript $\Theta$. Set $q=\pi^{-1}$. The normalizations are
\[
 E_{g,n}(\dd)=\frac{J_{g,n}(\dd)}{V_{g,n}},\qquad
 J_{g,n}(\dd)=4^{|\dd|}[\tau_{d_1}\cdots\tau_{d_n}]_g^\Theta,
\]
and
\begin{equation}
\begin{aligned}
 E_{g,n}(\dd)&\sim\sum_{r\ge0}P_r(\dd;n)g^{-r},\qquad P_0=1,\\
 a_n(g)&=\frac{\pi}{2}(2g-2+n)\frac{V_{g,n}}{V_{g,n+1}}
        \sim1+\sum_{r\ge1}a_n^r g^{-r},\\
 b_n(g)&=4\frac{V_{g-1,n+2}}{V_{g,n}}
        \sim1+\sum_{r\ge1}b_n^r g^{-r}.
\end{aligned}
\label{eq:coefficient-normalizations}
\end{equation}
The first author proved the existence of these fixed-marking asymptotic
expansions to arbitrary order \cite{Huang}. We extract their coefficients from the
exact recursion, rather than reuse its printed second-order
simplifications. All equalities between coefficient polynomials below
are exact over $\mathbb Q[n,q]$.

\subsection{The finite coefficient algorithm}

Put
\[
 \Delta_L=[s^L]\frac{1-s}{\cos(\pi\sqrt{s}/2)},\qquad
 \omega_k=\frac{(-1)^k(\pi/2)^{2k+1}}{(2k+1)!}.
\]
Subtracting the recursion for $J(\dd+\mathbf e_1)$ from that
for $J(\dd)$ gives
\begin{equation}
\begin{aligned}
 &J_{g,n}(\dd)-J_{g,n}(\dd+\mathbf e_1)\\
 &\quad=\sum_{j=2}^n(2d_j+1)\sum_{L\ge0}\Delta_L
 J_{g,n-1}(d_1+d_j+L,\dd_{\widehat1,\widehat j})\\
 &\qquad\quad+8\sum_{L\ge0}\Delta_L
 \sum_{u+v=L+d_1-1}\mathcal R_{g,n}(u,v;\dd_{\widehat1}).
\end{aligned}
\label{eq:coefficient-defect}
\end{equation}
The sum is finite at any fixed genus because excess degrees give zero.
The passage to coefficient moments is justified by the exponential
decay of $\Delta_L$ and the polynomial remainder bounds in
Theorem~4.1 of \cite{Huang}; see also
Appendix~\ref{sec:kdv-verification}. In particular, exponentially small
upper-tail errors do not affect any power of $g^{-1}$.

The two moment functionals required by the algorithm are
\begin{equation}
\begin{aligned}
 \mu_j&=\sum_{L\ge0}L^j\Delta_L
 =\left.(s\partial_s)^j
 \frac{1-s}{\cos(\pi\sqrt{s}/2)}\right|_{s=1},\\
 \nu_j&=\sum_{k\ge0}k^j\omega_k
 =\left.\frac\pi2\left(\frac z2\partial_z\right)^j
 \frac{\sin z}{z}\right|_{z=\pi/2}.
\end{aligned}
\label{eq:coefficient-moments}
\end{equation}
The apparent singularity in the first line is removable. Some useful
values, which also fix all normalization factors, are
\begin{equation}
\begin{gathered}
 \mu_0=4q,\qquad \mu_1=q,\qquad
 \mu_2=\frac q2+\frac1{12q},\\
 \nu_0=1,\qquad\nu_1=-\frac12,\qquad
 \nu_2=\frac14-\frac1{16q^2},\\
 \nu_3=-\frac18,\qquad
 \nu_4=\frac1{16}-\frac1{64q^2}+\frac1{256q^4}.
\end{gathered}
\label{eq:coefficient-first-moments}
\end{equation}
For a polynomial $f$, evaluate the infinite sum by replacing each
$L^j$ by $\mu_j$, or each $k^j$ by $\nu_j$. Threshold exceptions at
finitely many small integers are corrected before this substitution.

\paragraph{Volume ratios and genus shifts.}
Every volume ratio that appears in the defect is obtained from $a_n$
and $b_n$. Explicitly, with $m=2h+r\ge0$,
\begin{equation}
\begin{aligned}
 \frac{V_{g-h,n-r}}{V_{g,n}}
 ={}&4^{-h}\prod_{j=0}^{h-1}b_{n+2j}(g-j)\\
 &\times\prod_{k=0}^{m-1}
 \frac{a_{n-r+k}(g-h)}
 {(\pi/2)(2g-2h+n-r+k-2)}.
\end{aligned}
\label{eq:coefficient-child-ratio}
\end{equation}
An empty product is one. Here $r=-1$ is allowed when $h=1$, giving
the nonseparating child ratio. The leading order is $g^{-m}$.
For a known series $c(g)=1+\sum_{j\ge1}c_jg^{-j}$, genus shifting
uses the exact identity
\[
 [g^{-k}]c(g-h)
 =\sum_{j=1}^k\binom{k-1}{j-1}h^{k-j}c_j.
\]
Consequently the order-$R$ defect uses only ratios and descendant
coefficients of orders strictly below $R$. This is the triangular
feature that makes the computation finite.

\paragraph{Separating terms through order five.}
An end with small genus $h$ and $r$ external markings begins at
$g^{-(2h+r)}$. Thus only $h=1,2$ enter through $g^{-5}$. The exact
small-genus values needed are
\begin{equation}
\begin{aligned}
 J_{1,r+1}(0^{r+1})&=\frac{r!}{8},\\
 J_{2,r+1}(0^{r+1})&=\frac{3\pi^2(r+3)!}{128},\\
 J_{2,r+1}(1,0^r)&=\frac{9(r+2)!}{16}.
\end{aligned}
\label{eq:coefficient-low-genus}
\end{equation}
All other positive-degree genus-one values and total-degree-at-least-two
genus-two values vanish. These values follow from the initial theta
intersections and the exact recursion; the independent KdV appendix
derives them again. The two ordered ends of a separating sum give a
factor two. The factorial convolution estimate discards the remaining
genera only after showing they are beyond the requested order.

\paragraph{Integrating the descendant defect.}
For the remainder of this section write $D=|\dd|$ and
$M_j=\#\{i:d_i\ge j\}$. This $D$ is the total descendant degree,
not the unused degree in the finite-genus recursion. Through the
computed order $R=5$, the state needed for $P_R$ is
$(D,n,M_1,\ldots,M_{R-1})$. Raising a distinguished degree $x$ by
one changes $D$ to $D+1$ and $M_{x+1}$ to $M_{x+1}+1$.
The coefficient of $g^{-R}$ in
\eqref{eq:coefficient-defect}, divided by $V_{g,n}$, therefore
specifies the finite differences of $P_R$.

For $x\ge R-1$, the computed defect is independent of $x$.
Integrate first in $D$ using
\[
 \sum_{a=0}^{D-1}a^j
 =\frac{\mathsf B_{j+1}(D)-\mathsf B_{j+1}(0)}{j+1},
\]
where $\mathsf B_j$ denotes the Bernoulli polynomial. Correct the
cases $x=0,\ldots,R-2$ by discrete integration in $M_{x+1}$.
The normalization $P_R(0^n;n)=0$ fixes the integration constant.
For merge terms, one groups the partner degrees into the finitely
many values below the current threshold and one remaining group.
For node terms, the generic sum over $u+v=L+x-1$ is polynomial in
$L$; one subtracts the finitely many low-degree exceptions. This
reduces every sum to \eqref{eq:coefficient-moments}. All these
finite-difference identities have been checked exactly, including
the threshold cases, through order five.

\paragraph{Recovering the two ratio series.}
After $P_R$ is known, the dilaton and first KdV identities give
\begin{equation}
 a_n^R=\sum_{k\ge0}\omega_k
 P_R(k,0^n;n+1),
\label{eq:coefficient-marking-extraction}
\end{equation}
and
\begin{equation}
\begin{aligned}
 b_n^R={}&P_R(1,0^{n-1};n)\\
 &-24[g^{-R}]\sum_{h=1}^{g-1}\sum_{r=0}^{n-2}
 \binom{n-2}{r}\frac{V_{h,r+2}V_{g-h,n-r}}{V_{g,n}}.
\end{aligned}
\label{eq:coefficient-genus-extraction}
\end{equation}
Thus one proceeds in the order
\[
 (P_1,a^1,b^1),\ (P_2,a^2,b^2),\ \ldots,\ (P_5,a^5,b^5).
\]
The fifth ratio coefficients are necessary: fourth-order ratios
alone do not determine the genus-dependent constant in $d_n^4$.
The exact conversion from these five ratio orders to $d_n^1,\ldots,d_n^4$
is given below and proved generally in Section~\ref{sec:coefficient-structure}.

\subsection{The separating contribution and the printed discrepancy}

At order $g^{-2}$, the separating part of the adjacent defect has
only a genus-one component with no external marking. Its descendant
degree is zero. Both ordered ends contribute, and
\[
 J_{1,1}(0)=\frac18,\qquad
 \frac{V_{g-1,n}}{V_{g,n}}
 =\frac1{4\pi^2g^2}+O(g^{-3}).
\]
For distinguished degree $x$, the separating coefficient is therefore
\begin{equation}
\begin{aligned}
 &8\cdot2\cdot\frac18\cdot\frac1{4\pi^2}
 \sum_{L\ge0}\Delta_L\mathbf1_{\{L+x\ge1\}}\\
 &\qquad=\frac2{\pi^3}-\frac{\mathbf1_{\{x=0\}}}{2\pi^2}.
\end{aligned}
\label{eq:coefficient-separating-correction}
\end{equation}
Here $\Delta_0=1$ and $\sum_L\Delta_L=4/\pi$. \cite[equation~(3.23)]{Huang} displays one quarter of
\eqref{eq:coefficient-separating-correction}. There are also
inconsistencies in its subsequent second-order simplifications.
This is why simply iterating the printed $d_n^1$ would propagate an
incorrect normalization to the later coefficients. The calculation
below starts from \eqref{eq:coefficient-defect} instead.

\subsection{A worked calculation of the first correction}

The first two integrated descendant polynomials are
\begin{equation}
 P_1=q^2(-4D^2-4Dn+6D)
\label{eq:coefficient-P1}
\end{equation}
and
\begin{equation}
\begin{aligned}
 P_2={}&q^2\left(2D^2n-6D^2+2Dn^2-9Dn+9D+\frac{M_1}{2}\right)\\
 &-2Dq^3\\
 &+q^4\bigl(8D^4+16D^3n-40D^3+8D^2n^2-48D^2n\\
 &\hspace{22mm}+62D^2-12Dn^2+44Dn-38D\bigr).
\end{aligned}
\label{eq:coefficient-P2}
\end{equation}
For one positive degree $k$, the threshold variable is
$M_1=\mathbf1_{\{k\ge1\}}$. In particular the marking sum
of the $M_1/(2\pi^2)$ term is
\[
 \frac{q^2}{2}\sum_{k\ge1}\omega_k
 =\frac{q^2}{2}\left(1-\frac1{2q}\right).
\]
This is a finite threshold correction, and must not be replaced
by the generic value $M_1=1$ at $k=0$.
Substitution of \eqref{eq:coefficient-first-moments} into
\eqref{eq:coefficient-marking-extraction} and
\eqref{eq:coefficient-genus-extraction} gives
\begin{equation}
\begin{aligned}
 a_n^1={}&\frac14+(2n-2)q^2,\qquad b_n^1=(2-4n)q^2,\\
 a_n^2={}&-\frac n8+\frac9{32}-\frac q4
 +\left(-\frac32n^2+5n-3\right)q^2\\
 &+q^3+(8n^2-20n+12)q^4,\\
 b_n^2={}&(2n^2-7n+2)q^2-2q^3+(-4n^2+12n-8)q^4.
\end{aligned}
\label{eq:coefficient-first-ratios}
\end{equation}
Let $A_r(n)=[g^{-r}]\log a_n(g)$ and
$B_r(n)=[g^{-r}]\log b_n(g)$. At second order,
$A_2=a_n^2-(a_n^1)^2/2$ and $B_2=b_n^2-(b_n^1)^2/2$, hence
\begin{equation}
\begin{aligned}
 A_2(n)={}&-\frac n8+\frac14-\frac q4
 +\left(-\frac32n^2+\frac92n-\frac52\right)q^2\\
 &+q^3+(6n^2-16n+10)q^4,\\
 B_2(n)={}&(2n^2-7n+2)q^2-2q^3+(-12n^2+20n-10)q^4.
\end{aligned}
\label{eq:coefficient-first-log-ratios}
\end{equation}
The combination
\[
\begin{aligned}
 R_2(n)&=B_2(n-2)+A_2(n-2)+A_2(n-1)\\
 &=-\frac n4+\frac78-\frac q2-(n-1)(n-2)q^2
\end{aligned}
\]
enters the genus-shift identity, which gives
\begin{equation}
d_n^1=R_2(n)-\frac14
 =-\frac n4+\frac58-\frac q2-(n-1)(n-2)q^2.
\label{eq:coefficient-corrected-d1}
\end{equation}

\subsection{Detailed conversion at orders two, three, and four}

Define
\[
 R_j(n)=B_j(n-2)+A_j(n-2)+A_j(n-1).
\]
Running the finite coefficient algorithm to order five determines
$R_2,R_3,R_4,R_5$. To display the subsequent arithmetic without
several pages of redundant intermediate ratio polynomials, use
\begin{equation}
\begin{aligned}
 \ell_1&=R_2-\frac14,\\
 \ell_2&=\frac12\left(R_3-\frac16-\ell_1\right),\\
 \ell_3&=\frac13\left(R_4-\frac18-\ell_1-3\ell_2\right),\\
 \ell_4&=\frac14\left(R_5-\frac1{10}-\ell_1-4\ell_2-6\ell_3\right).
\end{aligned}
\label{eq:coefficient-four-log-steps}
\end{equation}
For example, the order-three output is equivalently
\[
 R_3=\frac16+\ell_1+2\ell_2,
\]
where
\[
\begin{aligned}
 \ell_2(n)={}&\frac{(n-3)(n-2)}{16}+\frac{2n-7}{8}q\\
 &+\frac{(n-2)(n^2-4n+2)}2q^2+(2-n)q^3\\
 &-(n-2)(n-1)(2n-5)q^4.
\end{aligned}
\]
The complete $\ell_3,\ell_4$ polynomials are supplied after the
explicit $d$ coefficients below. Finally exponentiate:
\begin{equation}
\begin{aligned}
 d_n^2&=\ell_2+\frac12\ell_1^2,\\
 d_n^3&=\ell_3+\ell_1\ell_2+\frac16\ell_1^3,\\
 d_n^4&=\ell_4+\ell_1\ell_3+\frac12\ell_2^2
       +\frac12\ell_1^2\ell_2+\frac1{24}\ell_1^4.
\end{aligned}
\label{eq:coefficient-exponentiation}
\end{equation}
These finite additions and multiplications yield the following
formulas. They were also computed independently from the one-positive-
descendant KdV method in Appendix~\ref{sec:kdv-verification}, which
checks the same coefficients without the multidescendant state variables.

\subsection{The four coefficients from the exact recursion}
Put $q=\pi^{-1}$. The following formulas use the normalization in Theorem~4.5, with its printed first coefficient corrected as explained above.
\begin{equation}\label{d1}
\begin{aligned}
d_n^{1}={}
&- \frac{2 n - 5}{8}
\\ &{}- \frac{q}{2}
\\ &{}- q^{2} \left(n - 2\right) \left(n - 1\right)
\end{aligned}
\end{equation}
\begin{equation}\label{d2}
\begin{aligned}
d_n^{2}={}
&\frac{12 n^{2} - 60 n + 73}{128}
\\ &{}+\frac{q \left(6 n - 19\right)}{16}
\\ &{}+\frac{q^{2} \left(6 n^{3} - 35 n^{2} + 59 n - 25\right)}{8}
\\ &{}+\frac{q^{3} \left(n - 3\right) \left(n - 2\right)}{2}
\\ &{}+\frac{q^{4} \left(n - 4\right) \left(n - 3\right) \left(n - 2\right) \left(n - 1\right)}{2}
\end{aligned}
\end{equation}
\begin{equation}\label{d3}
\begin{aligned}
d_n^{3}={}
&- \frac{5 \left(2 n - 5\right) \left(4 n^{2} - 20 n + 23\right)}{1024}
\\ &{}- \frac{q \left(60 n^{2} - 380 n + 637\right)}{256}
\\ &{}- \frac{q^{2} \left(180 n^{4} - 1560 n^{3} + 4591 n^{2} - 4845 n + 824\right)}{384}
\\ &{}- \frac{q^{3} \left(10 n^{3} - 77 n^{2} + 191 n - 153\right)}{16}
\\ &{}- \frac{q^{4} \left(n - 4\right) \left(n - 3\right) \left(30 n^{3} - 185 n^{2} + 321 n - 136\right)}{48}
\\ &{}- \frac{q^{5} \left(n - 5\right) \left(n - 4\right) \left(n - 3\right) \left(n - 2\right)}{4}
\\ &{}- \frac{q^{6} \left(n - 6\right) \left(n - 5\right) \left(n - 4\right) \left(n - 3\right) \left(n - 2\right) \left(n - 1\right)}{6}
\end{aligned}
\end{equation}
\begin{equation}\label{d4}
\begin{aligned}
d_n^{4}={}
&\frac{7 \left(80 n^{4} - 800 n^{3} + 2920 n^{2} - 4600 n + 2637\right)}{32768}
\\ &{}+\frac{q \left(280 n^{3} - 2660 n^{2} + 8918 n - 10981\right)}{2048}
\\ &{}+\frac{q^{2}}{3072}\bigl(840 n^{5} - 9660 n^{4} + 41174 n^{3}
\\ &\hspace{25mm}{}-75059 n^{2} + 39847 n + 20297\bigr)
\\ &{}+\frac{q^{3} \left(420 n^{4} - 4368 n^{3} + 15931 n^{2} - 23471 n + 11256\right)}{768}
\\ &{}+\frac{q^{4}}{768}\bigl(420 n^{6} - 6860 n^{5} + 44751 n^{4}
\\ &\hspace{25mm}{}-146762 n^{3} + 247173 n^{2} - 187490 n + 36330\bigr)
\\ &{}+\frac{q^{5} \left(n - 5\right) \left(n - 4\right) \left(42 n^{3} - 323 n^{2} + 793 n - 624\right)}{96}
\\ &{}+\frac{q^{6}}{48}\left(14 n^{3} - 91 n^{2} + 163 n - 71\right)
\prod_{j=3}^{6}(n-j)
\\ &{}+\frac{q^{7} \left(n - 7\right) \left(n - 6\right) \left(n - 5\right) \left(n - 4\right) \left(n - 3\right) \left(n - 2\right)}{12}
\\ &{}+\frac{q^{8}}{24}\prod_{j=1}^{8}(n-j)
\end{aligned}
\end{equation}
\subsection{Logarithmic coefficients}
For verification and generation, write $\log(1+\sum_{r\ge1}d_n^r x^r)=\sum_{r\ge1}\ell_r(n)x^r$.
\begin{equation}
\begin{aligned}
\ell_{1}(n)={}
&- \frac{2 n - 5}{8}
\\ &{}- \frac{q}{2}
\\ &{}- q^{2} \left(n - 2\right) \left(n - 1\right)
\end{aligned}
\end{equation}
\begin{equation}
\begin{aligned}
\ell_{2}(n)={}
&\frac{\left(n - 3\right) \left(n - 2\right)}{16}
\\ &{}+\frac{q \left(2 n - 7\right)}{8}
\\ &{}+\frac{q^{2} \left(n - 2\right) \left(n^{2} - 4 n + 2\right)}{2}
\\ &{}+q^{3} \left(2 - n\right)
\\ &{}- q^{4} \left(n - 2\right) \left(n - 1\right) \left(2 n - 5\right)
\end{aligned}
\end{equation}
\begin{equation}
\begin{aligned}
\ell_{3}(n)={}
&- \frac{\left(2 n - 5\right) \left(2 n^{2} - 10 n + 11\right)}{192}
\\ &{}- \frac{q \left(4 n^{2} - 28 n + 53\right)}{32}
\\ &{}- \frac{q^{2} \left(12 n^{4} - 108 n^{3} + 320 n^{2} - 312 n + 13\right)}{48}
\\ &{}+\frac{q^{3} \left(18 n^{2} - 96 n + 119\right)}{24}
\\ &{}+\frac{q^{4} \left(26 n^{4} - 224 n^{3} + 670 n^{2} - 796 n + 285\right)}{12}
\\ &{}- 2 q^{5} \left(n - 2\right) \left(2 n - 5\right)
\\ &{}- \frac{4 q^{6} \left(n - 2\right) \left(n - 1\right) \left(5 n^{2} - 27 n + 37\right)}{3}
\end{aligned}
\end{equation}
\begin{equation}
\begin{aligned}
\ell_{4}(n)={}
&\frac{\left(n - 3\right) \left(n - 2\right) \left(n^{2} - 5 n + 5\right)}{128}
\\ &{}+\frac{q \left(4 n^{3} - 42 n^{2} + 159 n - 227\right)}{64}
\\ &{}+\frac{q^{2}}{64}\bigl(8 n^{5} - 96 n^{4} + 416 n^{3}
\\ &\hspace{23mm}{}-720 n^{2} + 182 n + 479\bigr)
\\ &{}- \frac{q^{3} \left(36 n^{3} - 234 n^{2} + 478 n - 293\right)}{96}
\\ &{}-\frac{q^{4}}{24}\bigl(42 n^{5} - 492 n^{4} + 2164 n^{3}
\\ &\hspace{23mm}{}-4312 n^{2} + 3563 n - 641\bigr)
\\ &{}+\frac{q^{5} \left(58 n^{3} - 444 n^{2} + 1097 n - 879\right)}{12}
\\ &{}+\frac{q^{6}}6\bigl(68 n^{5} - 792 n^{4} + 3548 n^{3}
\\ &\hspace{23mm}{}-7524 n^{2} + 7313 n - 2409\bigr)
\\ &{}- 4 q^{7} \left(n - 2\right) \left(5 n^{2} - 27 n + 37\right)
\\ &{}- 2 q^{8} \left(n - 2\right) \left(n - 1\right) \left(14 n^{3} - 121 n^{2} + 355 n - 353\right)
\end{aligned}
\end{equation}

\section{Structure of the asymptotic coefficients}
\label{sec:coefficient-structure}

This section distinguishes consequences of the all-order ratio
polynomiality from patterns verified only through order four.
Set \(q=\pi^{-1}\), \(x=g^{-1}\), and write
\[
\begin{gathered}
 K_{g,n}=\frac{\Gamma(2g+n-2)\pi^{2g+n}}{2^n\sqrt g},
 \qquad V^\Theta_{g,n}\sim C K_{g,n}D_n(x),\\
 D_n(x)=1+\sum_{r\geq1}d_n^r x^r.
\end{gathered}
\]
All identities between asymptotic series below are formal identities in
\(x\), interpreted coefficient by coefficient. No convergence of these
infinite series is asserted.

\subsection{Exact recovery from the ratio coefficients}

Let
\[
\begin{aligned}
 a_n(x)&=1+\sum_{r\geq1}a_n^r x^r
 =\frac{(\pi/2)(2g-2+n)V^\Theta_{g,n}}{V^\Theta_{g,n+1}},\\
 b_n(x)&=1+\sum_{r\geq1}b_n^r x^r
 =\frac{4V^\Theta_{g-1,n+2}}{V^\Theta_{g,n}}.
\end{aligned}
\]
The normalizing factors cancel exactly, giving
\begin{equation}
 a_n(x)=\frac{D_n(x)}{D_{n+1}(x)},\qquad
 b_n(x)=(1-x)^{-1/2}
 \frac{D_{n+2}(x/(1-x))}{D_n(x)}.
 \label{eq:structure-ratio-identities}
\end{equation}
Define
\[
 \log D_n(x)=\sum_{r\geq1}\ell_r(n)x^r,\qquad
 \log a_n(x)=\sum_{r\geq1}A_r(n)x^r,\qquad
 \log b_n(x)=\sum_{r\geq1}B_r(n)x^r.
\]
Taking logarithms in \eqref{eq:structure-ratio-identities} yields
\begin{equation}
 \ell_r(n+1)-\ell_r(n)=-A_r(n)
 \label{eq:structure-marking-difference}
\end{equation}
and
\begin{equation}
 B_r(n)=\frac{1}{2r}+\ell_r(n+2)-\ell_r(n)
 +\sum_{j=1}^{r-1}\binom{r-1}{j-1}\ell_j(n+2).
 \label{eq:structure-genus-difference}
\end{equation}
Here
\[
 [x^r]\left(\frac{x}{1-x}\right)^j=\binom{r-1}{j-1}
 \qquad(1\leq j\leq r).
\]

Put \(R_j(n)=B_j(n-2)+A_j(n-2)+A_j(n-1)\).
Apply \eqref{eq:structure-genus-difference} at order \(r+1\) and at
\(n-2\), then use \eqref{eq:structure-marking-difference} to eliminate
\(\ell_{r+1}(n)-\ell_{r+1}(n-2)\). The resulting triangular formula is
\begin{equation}
\displaystyle
 \ell_r(n)=\frac1r\left[
 R_{r+1}(n)-\frac{1}{2(r+1)}
 -\sum_{j=1}^{r-1}\binom r{j-1}\ell_j(n)\right].
 \label{eq:structure-recovery}
\end{equation}
Initially this identity is obtained for integers \(n\geq2\).
The polynomial identities then extend it to \(n=0,1\).
In this extension, negative arguments of \(A_j\) or \(B_j\) mean
evaluation of their polynomials; they do not designate moduli spaces
with a negative number of markings.

For the first four orders,
\begin{equation}
\begin{aligned}
 \ell_1&=R_2-\frac14,\\
 \ell_2&=\frac12\left(R_3-\frac16-\ell_1\right),\\
 \ell_3&=\frac13\left(R_4-\frac18-\ell_1-3\ell_2\right),\\
 \ell_4&=\frac14\left(R_5-\frac1{10}-\ell_1-4\ell_2-6\ell_3\right).
\end{aligned}
\label{eq:structure-first-four-logs}
\end{equation}
Thus the ratio coefficients through \(g^{-5}\) suffice to determine
\(d_n^1,\ldots,d_n^4\). Truncation at \(g^{-4}\) leaves the
genus-normalization information needed for the fourth coefficient
undetermined.

Exponentiation gives the complete Bell-polynomial formula
\begin{equation}
 d_n^r=
 \sum_{\substack{m_1,m_2,\ldots\geq0\\\sum_{j\geq1}j m_j=r}}
 \prod_{j\geq1}\frac{\ell_j(n)^{m_j}}{m_j!},
 \qquad
 r d_n^r=\sum_{j=1}^r j\ell_j(n)d_n^{r-j},
 \quad d_n^0=1.
 \label{eq:structure-bell}
\end{equation}
In particular,
\begin{equation}
\begin{aligned}
 d_n^1&=\ell_1,\\
 d_n^2&=\ell_2+\frac12\ell_1^2,\\
 d_n^3&=\ell_3+\ell_1\ell_2+\frac16\ell_1^3,\\
 d_n^4&=\ell_4+\ell_1\ell_3+\frac12\ell_2^2
       +\frac12\ell_1^2\ell_2+\frac1{24}\ell_1^4.
\end{aligned}
\label{eq:structure-first-four-bell}
\end{equation}

\subsection{A general polynomial structure theorem}

The structural input is the ratio polynomiality established in
Theorem~4.1(2) of \cite{Huang}:
\(a_n^r,b_n^r\in\mathbb Q[n,q]\), with degree at most \(r\)
in \(n\) and at most \(2r\) in \(q\). The symbolic coefficient
computations here use the corrected values obtained from the exact
recursions. We now restate Theorem \ref{Structurethm} and prove it.

\begin{theorem}[Theorem 3]
For every \(r\geq1\),
\[
 \ell_r(n),d_n^r\in\mathbb Q[n,q],\qquad
 \deg_n\ell_r\leq r+1,\qquad \deg_n d_n^r=2r.
\]
The two highest coefficients in \(n\) are
\begin{equation}
\begin{aligned}
 \relax[n^{2r}]d_n^r&=\frac{(-1)^r q^{2r}}{r!},\\
 [n^{2r-1}]d_n^r&=
 \frac{(-1)^r q^{2r-2}}{(r-1)!}
 \left(\frac{2r-1}{4}-(2r+1)q^2\right).
\end{aligned}
 \label{eq:structure-top-two}
\end{equation}
The following bounds on inverse-pi degree hold:
\begin{equation}
 \deg_q\ell_r\leq2r+2,\qquad
 \deg_q\bigl([n^j]\ell_r\bigr)\leq2r\quad(j\geq1),\qquad
 \deg_q d_n^r\leq2r+2\lfloor r/2\rfloor.
 \label{eq:structure-q-bounds}
\end{equation}
These inverse-pi bounds are not claimed to be optimal.
\end{theorem}

\begin{proof}
A product contributing to \(A_r=[x^r]\log a_n(x)\) or
\(B_r=[x^r]\log b_n(x)\) has total asymptotic order \(r\).
The ratio polynomiality therefore implies
\[
 A_r,B_r\in\mathbb Q[n,q],\qquad
 \deg_n A_r,\deg_n B_r\leq r,\qquad
 \deg_q A_r,\deg_q B_r\leq2r.
\]
Equation~\eqref{eq:structure-recovery} proves polynomiality of
\(\ell_r\), inductively in \(r\). The forward difference of a
nonconstant polynomial of degree \(m\) has degree \(m-1\).
Equation~\eqref{eq:structure-marking-difference} consequently gives
\(\deg_n\ell_r\leq r+1\).

The explicitly computed first two logarithmic coefficients are
\begin{equation}
 \ell_1(n)=-\frac n4+\frac58-\frac q2-(n-1)(n-2)q^2
 \label{eq:structure-ell-one}
\end{equation}
and
\begin{equation}
\begin{aligned}
 \ell_2(n)={}&\frac{n^2}{16}-\frac{5n}{16}+\frac38
 +\left(\frac n4-\frac78\right)q\\
 &+\left(\frac{n^3}{2}-3n^2+5n-2\right)q^2
 +(2-n)q^3\\
 &+(-2n^3+11n^2-19n+10)q^4.
\end{aligned}
\label{eq:structure-ell-two}
\end{equation}
In a monomial of \eqref{eq:structure-bell}, its degree in \(n\) is at
most
\[
 \sum_{j\geq1}(j+1)m_j
 =r+\sum_{j\geq1}m_j\leq2r.
\]
Equality is possible only when \(m_1=r\). Since
\([n^2]\ell_1=-q^2\), this proves both the first formula in
\eqref{eq:structure-top-two} and \(\deg_n d_n^r=2r\).

For degree \(2r-1\), only \(\ell_1^r/r!\) and, when \(r\geq2\),
\(\ell_2\ell_1^{r-2}/(r-2)!\) can contribute. All other Bell
monomials have \(\sum m_j\leq r-2\). From
\eqref{eq:structure-ell-one}--\eqref{eq:structure-ell-two},
\[
 [n]\ell_1=-\frac14+3q^2,\qquad
 [n^3]\ell_2=\frac{q^2}{2}-2q^4.
\]
For \(r\geq2\), the desired coefficient is therefore
\[
 \frac{(-q^2)^{r-1}}{(r-1)!}
       \left(-\frac14+3q^2\right)
 +\frac{(-q^2)^{r-2}}{(r-2)!}
       \left(\frac{q^2}{2}-2q^4\right),
\]
which simplifies to the second formula in
\eqref{eq:structure-top-two}. The case \(r=1\) follows directly from
\eqref{eq:structure-ell-one}.

Equation~\eqref{eq:structure-recovery}, together with the bounds on
\(A_{r+1}\) and \(B_{r+1}\), gives \(\deg_q\ell_r\leq2r+2\).
For every power \(q^m\) with \(m>2r\), the coefficient of that power
in \eqref{eq:structure-marking-difference} has zero forward
difference. It is consequently constant in \(n\). This proves the
second bound in \eqref{eq:structure-q-bounds}.

Finally, \eqref{eq:structure-ell-one} gives the sharper bound
\(\deg_q\ell_1=2\). Each Bell monomial therefore has inverse-pi
degree at most
\[
 2m_1+\sum_{j\geq2}(2j+2)m_j
 =2r+2\sum_{j\geq2}m_j
 \leq2r+2\lfloor r/2\rfloor,
\]
proving the remaining claim.
\end{proof}

For example, \eqref{eq:structure-top-two} gives the exact leading
polynomial terms
\[
\begin{aligned}
 d_n^2&=\frac{q^4}{2}n^4
       +\left(\frac{3q^2}{4}-5q^4\right)n^3
       +P_2(n,q),\\
 d_n^3&=-\frac{q^6}{6}n^6
       +\left(-\frac{5q^4}{8}+\frac{7q^6}{2}\right)n^5
       +P_3(n,q),\\
 d_n^4&=\frac{q^8}{24}n^8
       +\left(\frac{7q^6}{24}-\frac{3q^8}{2}\right)n^7
       +P_4(n,q),
\end{aligned}
\]
where \(\deg_nP_r\leq2r-2\).
These statements concern polynomial coefficients. They do not extend
the fixed-marking genus asymptotic to growing \(n\) without additional
uniform estimates.

\subsection{Patterns verified through order four}

The exact computations at \(r=1,2,3,4\) have
\(\deg_q d_n^r=2r\) and satisfy
\begin{equation}
\begin{aligned}
 \relax[q^{2r}]d_n^r&=\frac{(-1)^r}{r!}
                  \prod_{j=1}^{2r}(n-j),\\
 [q^{2r-1}]d_n^r&=\frac{(-1)^r}{2(r-1)!}
                  \prod_{j=2}^{2r-1}(n-j).
\end{aligned}
\label{eq:structure-finite-q-patterns}
\end{equation}
The empty product is one. In this document,
\eqref{eq:structure-finite-q-patterns} and the sharper inverse-pi
degree \(2r\) are verified only for these four orders. An assertion for
every \(r\) would require a separate proof of the additional
highest-degree cancellations; it is not a consequence of
Theorem~\ref{Structurethm} alone.

\section{Large numbers of markings at fixed genus}
\label{sec:large-n}

In this section the genus is fixed and the number of markings tends
to infinity. We consider the cusp volumes \(V^\Theta_{g,n}\), as in
Theorem~\ref{superWP}. The theta degree condition makes this limit
particularly explicit: for \(g\geq2\), the quotient
\(V^\Theta_{g,n}/n!\) is exactly a polynomial in \(n\).
Consequently its expansion in inverse powers of \(n\) terminates.
No uniformity in a simultaneous limit of \(g\) and \(n\) is asserted.

\begin{theorem}[The complete expansion for large numbers of markings]
\label{large-n:theorem}
For every fixed integer \(g\geq2\), there are a positive rational
number \(\beta_g\) and rational numbers \(c_g^r\),
\(1\leq r\leq3g-4\), such that, for every integer \(n\geq1\),
\begin{equation}
 V^\Theta_{g,n}
 =(2\pi^2)^{g-1}\beta_g\,n!\,n^{3g-4}
 \left(1+\sum_{r=1}^{3g-4}\frac{c_g^r}{n^r}\right).
 \label{large-n:finite-expansion}
\end{equation}
In particular, for every fixed \(s\geq0\),
\begin{equation}
 V^\Theta_{g,n}
 =(2\pi^2)^{g-1}\beta_g\,n!\,n^{3g-4}
 \left(1+\sum_{r=1}^{s}\frac{c_g^r}{n^r}
       +O_{g,s}(n^{-s-1})\right),
 \qquad n\longrightarrow\infty,
 \label{large-n:asymptotic}
\end{equation}
where \(c_g^r=0\) for \(r>3g-4\). The remainder is identically
zero when \(s\geq3g-4\). The leading coefficient is
\begin{equation}
 \beta_g=
 \frac{\langle\tau_1^{g-1}\rangle_g^\Theta}
 {(g-1)!\,\Gamma(3g-3)}.
 \label{large-n:beta}
\end{equation}
All coefficients are given by the finite formula
\eqref{large-n:all-coefficients} below. In genera one and zero,
respectively, the exact formulas are
\[
 V^\Theta_{1,n}=\frac{(n-1)!}{8}\quad(n\geq1),
 \qquad
 V^\Theta_{0,n}=0\quad(n\geq3).
\]
\end{theorem}

\subsection{An exact formula from the theta dilaton identity}

Put \(a=2\pi^2\). A pure theta intersection will be denoted by
\[
 I_g(\boldsymbol m)
 =\left\langle\prod_{k=1}^{g-1}\tau_k^{m_k}\right\rangle_g^\Theta,
 \qquad
 |\boldsymbol m|=\sum_{k=1}^{g-1}m_k,
 \qquad
 \sum_{k=1}^{g-1}k m_k=g-1.
\]
For \(0\leq j\leq g-2\), define
\begin{equation}
 \rho_{g,j}=
 \frac{(-1)^j}{\Gamma(3g-3-j)}
 \sum_{\substack{m_1,\ldots,m_{g-1}\geq0\\
                  \sum k m_k=g-1\\
                  |\boldsymbol m|=g-1-j}}
 \frac{I_g(\boldsymbol m)}
      {\prod_{k=1}^{g-1}m_k!\,(k!)^{m_k}}.
 \label{large-n:rho}
\end{equation}
Set \(\rho_{g,j}=0\) outside this range. In particular,
\(\rho_{g,0}=\beta_g\). All the sums in this subsection are finite.

\begin{proposition}[Exact dependence on the marking number]
\label{large-n:gamma-proposition}
For \(g\geq2\) and every \(n\geq0\),
\begin{equation}
 V^\Theta_{g,n}
 =a^{g-1}\sum_{j=0}^{g-2}
       \rho_{g,j}\,\Gamma(n+3g-3-j).
 \label{large-n:gamma-sum}
\end{equation}
Equivalently, its exponential generating function is the rational
function
\begin{equation}
 \sum_{n\geq0}V^\Theta_{g,n}\frac{z^n}{n!}
 =a^{g-1}\sum_{j=0}^{g-2}
 \frac{\rho_{g,j}\Gamma(3g-3-j)}{(1-z)^{3g-3-j}}.
 \label{large-n:egf}
\end{equation}
The identity holds as a formal series and analytically for \(|z|<1\).
\end{proposition}

\begin{proof}
The exact kappa-to-descendant transformation in
Section~\ref{sec:interpolation} specializes to
\[
 V^\Theta_{g,n}
 =\left\langle\tau_0^n
  \exp\left(\sum_{k\geq1}\frac{(-1)^{k+1}a^k}{k!}\tau_k\right)
  \right\rangle_g^\Theta.
\]
The dimension condition is \(\sum k m_k=g-1\), independently of
the number of zero insertions. For \(g\geq2\), it forces
\(1\leq m:=|\boldsymbol m|\leq g-1\). Expanding the exponential
therefore gives the exact finite sum
\begin{equation}
 V^\Theta_{g,n}=a^{g-1}
 \sum_{\sum k m_k=g-1}
 \frac{(-1)^{g-1+m}}
      {\prod_k m_k!\,(k!)^{m_k}}
 \left\langle\tau_0^n\prod_k\tau_k^{m_k}\right\rangle_g^\Theta.
 \label{large-n:shift-sum}
\end{equation}
For \(N\) existing insertions, the theta dilaton identity reads
\[
 \left\langle\tau_0\prod_{i=1}^{N}\tau_{d_i}\right\rangle_g^\Theta
 =(2g-2+N)
 \left\langle\prod_{i=1}^{N}\tau_{d_i}\right\rangle_g^\Theta.
\]
Applying it \(n\) times yields
\begin{equation}
 \left\langle\tau_0^n\prod_k\tau_k^{m_k}\right\rangle_g^\Theta
 =\frac{\Gamma(n+2g-2+m)}{\Gamma(2g-2+m)}
   I_g(\boldsymbol m).
 \label{large-n:zero-insertions}
\end{equation}
All pairs involved are stable, including the pair with \(m\) markings
before the zero insertions. Substitute this identity into
\eqref{large-n:shift-sum} and put \(j=g-1-m\). Then
\((-1)^{g-1+m}=(-1)^j\) and \(2g-2+m=3g-3-j\), proving
\eqref{large-n:gamma-sum}. Finally, for a positive integer \(A\),
\[
 \sum_{n\geq0}\frac{\Gamma(n+A)}{n!}z^n
 =\Gamma(A)(1-z)^{-A}.
\]
Applying this to the finite gamma sum proves \eqref{large-n:egf}.
\end{proof}

\subsection{Every coefficient and the first two corrections}

For integers \(r,M\geq0\), let
\[
 E_r(M)=e_r(1,2,\ldots,M),\qquad E_0(M)=1,
\]
where \(e_r\) is the elementary symmetric polynomial; set
\(E_r(M)=0\) if \(r<0\) or \(r>M\). These rational integers can
be generated without any asymptotic calculation by
\[
 E_r(M)=E_r(M-1)+M E_{r-1}(M-1),\qquad E_0(0)=1.
\]
Write \(d=3g-4\). The complete coefficient formula is
\begin{equation}
 \displaystyle
 c_g^r=\frac1{\beta_g}
 \sum_{j=0}^{\min(r,g-2)}\rho_{g,j}E_{r-j}(d-j),
 \qquad 0\leq r\leq d,
 \label{large-n:all-coefficients}
\end{equation}
with \(c_g^0=1\). In particular,
\begin{equation}
\begin{aligned}
 c_g^1&=\frac{d(d+1)}2+\frac{\rho_{g,1}}{\beta_g},\\
 c_g^2&=\frac{d(d+1)(d-1)(3d+2)}{24}
       +\frac{d(d-1)}2\frac{\rho_{g,1}}{\beta_g}
       +\frac{\rho_{g,2}}{\beta_g}.
\end{aligned}
 \label{large-n:first-two}
\end{equation}
The first two partition sums beyond the leading term are
\begin{equation}
 \rho_{g,1}=
 -\frac{\langle\tau_2\tau_1^{g-3}\rangle_g^\Theta}
 {2(g-3)!\,\Gamma(3g-4)}\qquad(g\geq3)
 \label{large-n:rho-one}
\end{equation}
and
\begin{equation}
 \rho_{g,2}=
 \frac1{\Gamma(3g-5)}
 \left(
 \frac{\langle\tau_3\tau_1^{g-4}\rangle_g^\Theta}{6(g-4)!}
 +\frac{\langle\tau_2^2\tau_1^{g-5}\rangle_g^\Theta}{8(g-5)!}
 \right)\qquad(g\geq4),
 \label{large-n:rho-two}
\end{equation}
where the second summand is omitted when \(g=4\).
For smaller genera the corresponding \(\rho\)'s are zero.
The factors \(1/2\), \(1/6\), and \(1/8\) include both the
factorials in the shifts and the factorials for repeated insertions.

For completeness, the pure intersections in \eqref{large-n:rho}
can be computed using the unshifted super DVV recursion
\cite{NP1,HLX}. In the present conventions it is
\begin{equation}
\begin{aligned}
 &(2d_1+1)!!
   \left\langle\prod_{i=1}^{N}\tau_{d_i}\right\rangle_g^\Theta\\
 &\quad=\sum_{j=2}^{N}
 \frac{(2d_1+2d_j+1)!!}{(2d_j-1)!!}
 \left\langle\tau_{d_1+d_j}
       \prod_{i\ne1,j}\tau_{d_i}\right\rangle_g^\Theta\\
 &\qquad\quad+\frac12\sum_{u+v=d_1-1}(2u+1)!!(2v+1)!!
 \left[
 \left\langle\tau_u\tau_v\prod_{i=2}^{N}\tau_{d_i}
       \right\rangle_{g-1}^\Theta\right.\\
 &\hspace{49mm}\left.
 +\sum_{\substack{h+h'=g\\I\sqcup J=\{2,\ldots,N\}}}
 \left\langle\tau_u\prod_{i\in I}\tau_{d_i}\right\rangle_h^\Theta
 \left\langle\tau_v\prod_{i\in J}\tau_{d_i}\right\rangle_{h'}^\Theta
 \right].
\end{aligned}
 \label{large-n:pure-dvv}
\end{equation}
Here \(u,v\geq0\), splittings are ordered, \((-1)!!=1\), and
genus-zero and unstable terms are zero. The exceptional initial value
is \(\langle\tau_0\rangle_1^\Theta=1/8\); the recursion is used
for all other admissible pairs. Each term reduces \(2g-2+N\), so
this is a finite procedure. Equations~\eqref{large-n:rho} and
\eqref{large-n:all-coefficients} thus determine every coefficient
in \eqref{large-n:asymptotic} by exact rational arithmetic.

\begin{proof}[Proof of Theorem~\ref{large-n:theorem}]
For each term of \eqref{large-n:gamma-sum}, the integer shift gives
the exact identity
\[
 \frac{\Gamma(n+3g-3-j)}{n!}
 =\prod_{k=1}^{d-j}(n+k).
\]
Hence
\begin{equation}
 \frac{V^\Theta_{g,n}}{a^{g-1}n!}
 =\sum_{j=0}^{g-2}\rho_{g,j}
          \prod_{k=1}^{d-j}(n+k).
 \label{large-n:polynomial}
\end{equation}
The coefficient of \(n^d\) is \(\rho_{g,0}=\beta_g\), whose
strict positivity is proved in the next subsection. Expanding each
product and collecting \(n^{d-r}\) gives precisely
\eqref{large-n:all-coefficients}. This proves the exact polynomial
identity \eqref{large-n:finite-expansion}, its degree, and rationality
of all relative coefficients. The first two elementary symmetric sums
are
\[
 E_1(d)=\frac{d(d+1)}2,\qquad
 E_2(d)=\frac{d(d+1)(d-1)(3d+2)}{24},
\]
which prove \eqref{large-n:first-two}. Truncating a fixed polynomial
in \(n^{-1}\) proves \eqref{large-n:asymptotic} with the stated
remainder. For \(g=1\), the degree condition excludes every positive
descendant and every positive power of \(\kappa_1\). Repeated theta
dilaton, starting with \(1/8\) at \(n=1\), gives
\(V^\Theta_{1,n}=(n-1)!/8\). The genus-zero assertion follows from
\(\Theta_{0,n}=0\).
\end{proof}

\subsection{A positive recurrence for the leading constant}

The coefficient \(\beta_g\) in \eqref{large-n:beta} can be computed
without the full partition sum. Define \(b_1=1/8\) and
\(b_g=\Gamma(3g-1)\beta_g\) for \(g\geq2\).
Then
\begin{equation}
 \displaystyle
 (g-1)b_g=
 \frac{(3g-4)(3g-3)(3g-2)}{12}\,b_{g-1}
 +\frac{3g-2}{2}\sum_{h=1}^{g-1}b_hb_{g-h},
 \qquad g\geq2.
 \label{large-n:leading-recursion}
\end{equation}
Every coefficient on the right is positive. The recurrence proves
\(b_g>0\), and hence \(\beta_g>0\), by induction.

To derive it, retain only the two descendant times \(t_0=x\) and
\(t_1=s\) in the unshifted theta free energy \(\mathcal F\), and put
\(U=\hbar^2\partial_x^2\mathcal F\). This \(s\) is a descendant
time, unrelated to the truncation order in
\eqref{large-n:asymptotic}. At genus \(g\), the degree condition
forces exactly \(g-1\) insertions of \(\tau_1\). The theta dilaton
identity therefore gives
\begin{equation}
 U(x,s;\hbar)
 =\sum_{g\geq1}\hbar^{2g}b_g s^{g-1}
                   (1-x)^{-(3g-1)}.
 \label{large-n:two-time}
\end{equation}
For \(g\geq2\), its coefficient at \(x=0\) is
\[
 \frac{\langle\tau_0^2\tau_1^{g-1}\rangle_g^\Theta}{(g-1)!}
 =(3g-3)(3g-2)
   \frac{\langle\tau_1^{g-1}\rangle_g^\Theta}{(g-1)!}
 =b_g;
\]
for \(g=1\), it is \(\langle\tau_0^2\rangle_1^\Theta=1/8\).
The first KdV equation in these conventions is
\[
 \partial_s U=U\partial_xU+\frac{\hbar^2}{12}\partial_x^3U.
\]
The coefficient of
\(\hbar^{2g}s^{g-2}(1-x)^{-(3g-1)}\) on its left is
\((g-1)b_g\). On the right, the third derivative contributes the
first term of \eqref{large-n:leading-recursion}, and the quadratic
term contributes
\[
 \sum_{h=1}^{g-1}(3(g-h)-1)b_hb_{g-h}
 =\frac{3g-2}{2}\sum_{h=1}^{g-1}b_hb_{g-h}.
\]
This proves the recurrence with all normalization factors fixed.
For example,
\[
 \beta_2=\frac3{256},\qquad
 \beta_3=\frac{21}{40960},\qquad
 \beta_4=\frac{2407}{165150720},\qquad
 \beta_5=\frac{536219}{1743991603200}.
\]

\subsection{Low-genus examples}

For \(g=2\), there is only the partition \(1\), with
\(\langle\tau_1\rangle_2^\Theta=3/128\). Thus
\begin{equation}
 V^\Theta_{2,n}
 =\frac{3\pi^2}{128}(n+2)!
 =\frac{3\pi^2}{128}\,n!\,n^2
   \left(1+\frac3n+\frac2{n^2}\right).
 \label{large-n:genus-two}
\end{equation}
For \(g=3\), the two partitions are \(1+1\) and \(2\).
The pure recursion gives
\[
 \langle\tau_1^2\rangle_3^\Theta=\frac{63}{512},\qquad
 \langle\tau_2\rangle_3^\Theta=\frac{15}{1024}.
\]
Substitution into \eqref{large-n:rho} gives
\begin{equation}
 V^\Theta_{3,n}
 =(2\pi^2)^2\left(
   \frac{21}{40960}\Gamma(n+6)
  -\frac5{16384}\Gamma(n+5)\right).
 \label{large-n:genus-three-gamma}
\end{equation}
Equivalently, its full expansion is
\begin{equation}
\begin{aligned}
 V^\Theta_{3,n}
 ={}&(2\pi^2)^2\frac{21}{40960}\,n!\,n^5\left(
 1+\frac{605}{42n}+\frac{1660}{21n^2}
   +\frac{1225}{6n^3}\right.\\
 &\hspace{46mm}\left.
   +\frac{5129}{21n^4}+\frac{740}{7n^5}\right).
\end{aligned}
 \label{large-n:genus-three-expansion}
\end{equation}
For a further check, genus four gives
\begin{equation}
\begin{aligned}
 V^\Theta_{4,n}=(2\pi^2)^3\left(
  \frac{2407}{165150720}\Gamma(n+9)
 -\frac{575}{22020096}\Gamma(n+8)
 +\frac{35}{9437184}\Gamma(n+7)\right).
\end{aligned}
 \label{large-n:genus-four-gamma}
\end{equation}
These identities hold for all \(n\geq0\) in their gamma or
polynomial forms. The representations with explicit negative powers
of \(n\) are used for \(n\geq1\).

\appendix
\section{Independent verification by the KdV hierarchy}
\label{sec:kdv-verification}

This section gives a second exact procedure for computing the four
coefficients. It uses only one positive descendant, and is independent
of the multidescendant defect calculation. Put
\[
q=\pi^{-1},\qquad x=g^{-1},\qquad
P_{g,n}=C\frac{\Gamma(2g+n-2)\pi^{2g+n}}{2^n\sqrt g},
\]
and write the formal asymptotic series
\[
D_n(x)=\frac{V_{g,n}}{P_{g,n}}
=\sum_{r\geq0}d_r(n)x^r,\qquad d_0(n)=1,
\qquad
F_{n,k}(x)=\frac{J_{g,n}(k,0^{n-1})}{P_{g,n}}.
\]
Here \(d_r(n)=d_n^r\), and \(F_{n,0}(x)=D_n(x)\).
All computations take place over \(\mathbb Q[n,q]\); no floating-point
approximation or numerical recognition of constants is used.

\subsection{An exact recurrence with one positive descendant}

The KdV hierarchy, with the normalization
\(p_k=16^k(2k+1)!!\), gives, for \(k\geq1\) and \(n\geq3\),
\begin{equation}
\begin{aligned}
J_{g,n}(k,0^{n-1})
={}&4J_{g-1,n+2}(k-1,0^{n+1})\\
&+16\sum_{h=1}^{g-1}\sum_{r=0}^{n-3}\binom{n-3}{r}
 J_{h,r+2}(k-1,0^{r+1})V_{g-h,n-r}\\
&+32\sum_{h=1}^{g-1}\sum_{r=0}^{n-3}\binom{n-3}{r}
 J_{h,r+3}(k-1,0^{r+2})V_{g-h,n-r-1}.
\end{aligned}
\label{eq:kdv-one-positive}
\end{equation}
To check the normalization, the unnormalized Lenard relation is
\[
(2k+1)\langle\!\langle\tau_k\tau_0^2\rangle\!\rangle
=\langle\!\langle\tau_{k-1}\tau_0\rangle\!\rangle
 \langle\!\langle\tau_0^3\rangle\!\rangle
+2\langle\!\langle\tau_{k-1}\tau_0^2\rangle\!\rangle
 \langle\!\langle\tau_0^2\rangle\!\rangle
+\frac{\hbar^2}{4}
 \langle\!\langle\tau_{k-1}\tau_0^4\rangle\!\rangle.
\]
Since \(p_k/((2k+1)p_{k-1})=16\), its three coefficients become
\(16,32,4\). Differentiating another \(n-3\) times in the zero
descendant variable gives the binomial factors in
\eqref{eq:kdv-one-positive}. At \(k=1\), symmetry combines the two
product sums into the integrated KdV identity used in the main proof.

\subsection{The two exact time identities}

Let \(a=\pi/2\), and define
\[
c_k=\frac{(-1)^k a^{2k}}{(2k)!},\qquad
w_k=\frac{c_k}{2k+1},\qquad W_k=a w_k.
\]
The theta dilaton identity and the weighted BGW string identity give
\begin{align}
(2g-2+n)V_{g,n}
 &=\sum_{k\geq0}w_kJ_{g,n+1}(k,0^n),
\label{eq:kdv-time-one}\\
nV_{g,n}
 &=\sum_{k\geq0}c_kJ_{g,n+1}(k,0^n).
\label{eq:kdv-time-two}
\end{align}
The second equation also follows from the first one and homogeneity:
\(\sum kT_k\partial_{T_k}V_{g,n}=(g-1)V_{g,n}\).
Thus it does not introduce a new assumption about the volumes.
After dividing by \(P_{g,n+1}\), these become
\begin{align}
D_n(x)&=\sum_{k\geq0}W_kF_{n+1,k}(x),
\label{eq:kdv-normal-time-one}\\
\frac{nqx}{1+(n-2)x/2}D_n(x)
 &=\sum_{k\geq0}c_kF_{n+1,k}(x).
\label{eq:kdv-normal-time-two}
\end{align}
Both identities are needed to determine the polynomial dependence on
the number of markings efficiently.

\subsection{The finite amount of low-genus data required}

If one component has genus \(h\) and \(s\) markings, the other has
genus \(g-h\) and \(n+2-s\) markings. Put \(m=2h+s-2\). The ratio of
the normalizing factors is exactly
\begin{equation}
\frac{P_{g-h,n+2-s}}{P_{g,n}}
=\frac{q^m x^m}{2^{2h}}
 (1-hx)^{-1/2}
 \prod_{\ell=1}^{m}
 \left(1+\frac{n-2-\ell}{2}x\right)^{-1}.
\label{eq:kdv-small-factor}
\end{equation}
Consequently a component can contribute at order at most \(x^5\) only
when \(2h+s-2\leq5\). The only possibilities are
\[
h=1,\quad 2\leq s\leq5;
\qquad
h=2,\quad 2\leq s\leq3.
\]
They require just the following exact values:
\begin{equation}
\begin{aligned}
V_{1,s}&=\frac{(s-1)!}{8},
& J_{1,s}(k,0^{s-1})&=0\quad(k\geq1),\\
V_{2,s}&=\frac{3(s+2)!}{128q^2},
& J_{2,s}(1,0^{s-1})&=\frac{9(s+1)!}{16},\\
&&J_{2,s}(k,0^{s-1})&=0\quad(k\geq2).
\end{aligned}
\label{eq:kdv-low-genus}
\end{equation}
These follow from \(\langle\tau_0\rangle_1^\Theta=1/8\),
\(\langle\tau_1\rangle_2^\Theta=3/128\), and the theta dilaton
identity. In particular, none of the required low-genus normalization
is inferred from an asymptotic formula.

To justify truncation, fix the marking count \(s\) on the smaller
component. Equations~\eqref{eq:convolution} and
Lemma~\ref{uniform:factorial-convolution}, with an arbitrary fixed
lower cutoff \(H\), bound the sum over \(h\geq H\), \(h\leq g/2\),
after division by \(V_{g,n}\), by \(O(g^{2-s-2H})\).
Choose the least \(H\geq1\) with \(2H+s-2\geq6\).
Every omitted endpoint contribution is then \(O(g^{-6})\), as is
the whole middle range. The same argument applies after interchanging
the components. Descendant monotonicity makes this estimate uniform
in \(k\). Thus only the listed components can affect orders through five.

\subsection{Formal coefficient recursion}

Write \(F_{n,k}(x)=\sum_{r\geq0}F_r(n,k)x^r\). The nonseparating term
in \eqref{eq:kdv-one-positive} becomes
\[
(1-x)^{-1/2}F_{n+2,k-1}\left(\frac{x}{1-x}\right).
\]
Thus its contribution at order \(r\) is
\begin{equation}
\sum_{j=0}^{r}
\frac{(j+\tfrac12)_{r-j}}{(r-j)!}
F_j(n+2,k-1),
\label{eq:kdv-nonseparating-series}
\end{equation}
where \((z)_m\) denotes the rising factorial.

For the product terms, let
\[
\begin{aligned}
M_+(n,s)&=\binom{n-3}{s-2}+2\binom{n-3}{s-3},\\
M_0(n,s)&=\binom{n-3}{s-3}+2\binom{n-3}{s-2}.
\end{aligned}
\]
The first multiplicity applies when the positive descendant is on the
small-genus component; the second applies when the cusp component is
the small one. A binomial coefficient with a negative lower index is
zero. For \(m=2h+s-2\), define
\[
K_{r,m,j,h}(n)=\frac{q^m}{2^{2h}}
[x^{r-m-j}](1-hx)^{-j-1/2}
\prod_{\ell=1}^{m}
\left(1+\frac{n-2-\ell}{2}x\right)^{-1}.
\]
The coefficient of a negative power is defined to be zero.
Then the full product contribution to \(F_r(n,k)\) is the finite sum
\begin{equation}
\begin{aligned}
16\sum_{\substack{h\geq1,\ s\geq2\\m=2h+s-2\leq r}}
\sum_{j=0}^{r-m}K_{r,m,j,h}(n)
\bigl[{}&M_+(n,s)J_{h,s}(k-1,0^{s-1})d_j(n+2-s)\\
&+M_0(n,s)V_{h,s}F_j(n+2-s,k-1)\bigr].
\end{aligned}
\label{eq:kdv-product-series}
\end{equation}
Equations~\eqref{eq:kdv-nonseparating-series} and
\eqref{eq:kdv-product-series}, together with \(F_r(n,0)=d_r(n)\),
are the exact finite-order recursion implemented in the accompanying
script.

\subsection{Summing over the descendant degree}

For fixed order \(r\leq5\), \(F_r(n,k)\) is a polynomial in \(k\)
of degree at most \(2r\) for \(k\geq\lfloor r/2\rfloor\).
The smaller values are treated separately. The exceptional range
comes from the dimension bounds on
\(J_{h,s}(k-1,0^{s-1})\) in the small component.

The required moments are elementary. If
\(\mathcal D_a=(a/2)\partial_a\), then
\begin{equation}
\sum_{k\geq0}W_k k^j
=\left.a\mathcal D_a^j\left(\frac{\sin a}{a}\right)
 \right|_{a=\pi/2},\qquad
\sum_{k\geq0}c_k k^j
=\left.\mathcal D_a^j\cos a\right|_{a=\pi/2}.
\label{eq:kdv-moment-evaluation}
\end{equation}
In particular,
\[
\sum W_k=1,\quad\sum kW_k=-\frac12,
\qquad\sum c_k=0,\quad\sum kc_k=-\frac\pi4.
\]
If \(P_r(n,k)\) denotes the generic polynomial and \(k_0=\lfloor
r/2\rfloor\), each infinite sum is evaluated as its polynomial moment
sum plus the finite correction
\[
\sum_{0\leq k<k_0}\omega_k
\bigl(F_r(n,k)-P_r(n,k)\bigr),
\qquad \omega_k=W_k\ \text{or}\ c_k.
\]
These finite corrections are important: omitting them changes the odd
powers of \(\pi^{-1}\).

The contractions in \eqref{eq:kdv-normal-time-one}--
\eqref{eq:kdv-normal-time-two} are valid asymptotically as well.
The precise input is Theorem~4.1(1), equations~(4.1)--(4.2), of
\cite{Huang}. For fixed \(n,s\), its
polynomial remainder bound implies
\[
\left|\frac{J_{g,n}(k,0^{n-1})}{V_{g,n}}
-\sum_{j=0}^{s}\frac{e^j_{n,(k,0^{n-1})}}{g^j}\right|
\leq \frac{C_{n,s}(1+k)^{2s+2}}{g^{s+1}},
\qquad e^0_{n,\mathbf d}=1.
\]
Multiplication by the fixed-marking expansion of \(V_{g,n}/P_{g,n}\)
gives a polynomial remainder bound for \(F_{n,k}\) as well.
Both \(W_k\) and \(c_k\) decay factorially in absolute value, so
contracting this bound proves termwise summation with remainder
\(O_{n,s}(g^{-s-1})\). The proof of Theorem~4.27 uses the exact
recursion and the order-zero and order-one bounds; it does not require
the incorrectly simplified printed second coefficient.

The script constructs the generic polynomial from exact rational
values of \(k\), but does not rely on interpolation as an assumption.
It substitutes the resulting polynomial into
\eqref{eq:kdv-nonseparating-series}--\eqref{eq:kdv-product-series}
with \(k\) still symbolic and verifies that the difference vanishes
identically. The finitely many boundary values are computed directly.

\subsection{Why one extra order determines the constant term}

At order \(r\), the dependence of \(F_r(n,k)\) on the new coefficient
\(d_r\) is exactly \(d_r(n+2k)\); all other terms are already known
from lower orders. Put \(\partial=\partial_n\). On polynomials in
\(n\), the two coefficient equations involve the operators
\[
\mathcal A=1-\sin(ae^{\partial}),\qquad
\mathcal B=e^{\partial}\cos(ae^{\partial}),
\qquad a=\frac\pi2.
\]
These are formal Taylor series in \(\partial\); their action on each
polynomial is a finite sum. Their leading terms are
\[
\mathcal A=\frac{a^2}{2}\partial^2+O(\partial^3),\qquad
\mathcal B=-a\partial+O(\partial^2).
\]
Thus the weighted string equation determines the nonconstant part of
\(d_r(n)\), while the two equations impose a compatibility condition
on the constant term of \(d_{r-1}(n)\).

More explicitly, let \(c\) be the still undetermined constant term of
\(d_{r-1}\). Its contribution to \(F_r(n,k)\) is
\((r-\tfrac12)kc\). Its contributions to the right-hand sides of the
\(\mathcal A\)- and \(\mathcal B\)-equations are therefore
\[
-\frac{r-1/2}{2}c,
\qquad
qn c+\frac{a(r-1/2)}2c.
\]
Since \(\mathcal A/\mathcal B=-a\partial/2+O(\partial^2)\), the
coefficient of \(c\) in the compatibility condition is
\(-\frac12(r-1)\), which is nonzero for \(r\geq2\).
This determines \(c\) uniquely. Consequently order five is needed to
complete \(d_4\), even though no coefficient beyond \(d_4\) is being
reported.

These unique formal coefficients are the actual asymptotic
coefficients: Theorem~1.5 of \cite{Huang} supplies their existence
and polynomiality, while the remainder and convolution bounds just
proved show that they satisfy this finite system. One may first work
at each fixed integer \(n\geq20\), ensuring that every marking count
used through order five is admissible. Polynomial identity then extends
the resulting formulas to every \(n\geq0\).

For example, order one gives
\[
d_1(n)=-q^2n^2+\left(3q^2-\frac14\right)n+c_1.
\]
The order-two compatibility condition gives
\[
c_1=\frac58-\frac q2-2q^2,
\]
and hence the corrected first coefficient
\[
d_n^1=-\frac n4+\frac58-\frac1{2\pi}
-\frac{n^2-3n+2}{\pi^2}.
\]
The next three compatibility calculations give
\begin{align*}
c_2={}&12q^4+3q^3-\frac{25}{8}q^2-\frac{19}{16}q+\frac{73}{128},\\
c_3={}&-120q^6-30q^5+34q^4+\frac{153}{16}q^3
-\frac{103}{48}q^2-\frac{637}{256}q+\frac{575}{1024},\\
c_4={}&1680q^8+420q^7-\frac{1065}{2}q^6-130q^5
+\frac{6055}{128}q^4+\frac{469}{32}q^3\\
&+\frac{20297}{3072}q^2-\frac{10981}{2048}q
+\frac{18459}{32768},
\end{align*}
where \(c_j=d_j(0)\). The accompanying exact script provides all
nonconstant terms as well and verifies the recurrence through order
five.

\section{Independent exact and numerical checks}
\label{sec:numerical-checks}

The coefficient calculation can be checked against volumes obtained from
a different implementation of the full BGW string and KdV equations.
This section reports such checks. The symbolic procedure in the preceding
sections determines the coefficients; numerical extrapolation is used
only as an independent verification, not as a proof.

\subsection{A recurrence for the actual volumes}

Set the kappa parameter to \(a=1\), let \(x=t_0\), and put
\(y=(1-x)^{-1}\). For the shifted BGW free energy write
\[
 U=\hbar^2\partial_x^2\mathcal F
   =\sum_{g\geq1}\hbar^{2g}U_g(y),
 \qquad U_1(y)=\frac{y^2}{8},
 \qquad T_k=\frac{(-1)^{k+1}}{k!}.
\]
The Gelfand--Dikii polynomials satisfy
\[
 R_1=U,\qquad
 \partial_xR_{k+1}
 =\frac{U'R_k+2U R_k'+\hbar^2R_k'''/4}{2k+1}.
\]
Here \(\partial_xy^p=py^{p+1}\), so integration is performed
coefficient by coefficient in \(y\). The weighted BGW string identity
then gives
\[
 (1-x)U_g'-2U_g
 =\sum_{k=1}^{g-1}T_k
   [\hbar^{2g}](U'R_k+2U R_k'+\hbar^2R_k'''/4).
\]
Its right-hand side uses only lower genera. A coefficient of \(y^p\)
on the right is divided by \(p-2\) to obtain the corresponding
coefficient of \(U_g\). In particular,
\[
 U_2=\frac9{32}y^5,\qquad
 U_3=-\frac{225}{1024}y^7+\frac{1323}{512}y^8.
\]
Writing \(U_g=\sum_m u_{g,m}y^{2g+m}\), the physical volume is
\[
 V^\Theta_{g,n}
 =(2\pi^2)^{g-1}\sum_m u_{g,m}(2g+m)_{n-2},\qquad n\geq2,
\]
where \((z)_j=z(z+1)\cdots(z+j-1)\), and \((z)_0=1\).

As an exact low-genus check, the values of
\(V^\Theta_{g,2}/(2\pi^2)^{g-1}\) are
\[
\begin{array}{c|r}
g&V^\Theta_{g,2}/(2\pi^2)^{g-1}\\ \hline
2&9/32\\
3&2421/1024\\
4&89215/2048\\
5&1480147661/1048576\\
6&149143143477/2097152\\
7&690735643405217/134217728
\end{array}
\]
These were also calculated independently using the unshifted super DVV
recursion and the finite kappa-to-descendant partition sum. The two
calculations agree exactly, including all rational normalization factors.

\subsection{Checking four asymptotic coefficients}

Define the normalized volume
\[
 A_{g,n}:=
 \frac{2^n\sqrt g\,V^\Theta_{g,n}}
 {\pi^{-9/2}\Gamma(2g+n-2)\pi^{2g+n}}.
\]
The recurrence above was evaluated through genus \(85\) with
\(150\)-digit decimal arithmetic. To independently estimate the first
four coefficients, solve the linear system
\[
 A_{g,n}-1=\sum_{j=1}^{20}\widehat d_n^j g^{-j},
 \qquad g=47,49,\ldots,85.
\]
No symbolic coefficient is used in this fit. The following table
compares its results with the exact formulas evaluated at \(q=\pi^{-1}\).
The final column is the observed absolute discrepancy before rounding
the two displayed numerical values.

\begin{center}
\small
\begin{tabular}{@{}rr r r r@{}}
\hline
\(n\)&\(j\)&Exact formula (decimal)&Independent fit&Discrepancy\\
\hline
2&1&\(-0.03415494309190\)&\(-0.03415494309190\)&\(8.73\cdot10^{-20}\)\\
2&2&\(-0.11878292725012\)&\(-0.11878292725012\)&\(1.13\cdot10^{-16}\)\\
2&3&\(-0.12802915473069\)&\(-0.12802915473062\)&\(6.94\cdot10^{-14}\)\\
2&4&\(-0.20647548585120\)&\(-0.20647548587801\)&\(2.68\cdot10^{-11}\)\\
3&1&\(-0.48679731037657\)&\(-0.48679731037657\)&\(8.54\cdot10^{-20}\)\\
3&2&\(-0.02474701584178\)&\(-0.02474701584178\)&\(1.11\cdot10^{-16}\)\\
3&3&\(-0.05301804346794\)&\(-0.05301804346787\)&\(6.79\cdot10^{-14}\)\\
3&4&\(-0.16205647975330\)&\(-0.16205647977951\)&\(2.62\cdot10^{-11}\)\\
4&1&\(-1.14208204494592\)&\(-1.14208204494592\)&\(8.28\cdot10^{-20}\)\\
4&2&\(0.77031605230086\)&\(0.77031605230086\)&\(1.07\cdot10^{-16}\)\\
4&3&\(-0.53737666903133\)&\(-0.53737666903127\)&\(6.59\cdot10^{-14}\)\\
4&4&\(0.10807808654314\)&\(0.10807808651771\)&\(2.54\cdot10^{-11}\)\\
\hline
\end{tabular}
\end{center}

These discrepancies are not certified extrapolation error bounds.
They are an independent numerical agreement with the symbolic answer.
The decreasing accuracy for successive coefficients is expected:
extracting the coefficient of \(g^{-j}\) amplifies the effect of omitted
higher-order terms.

A complementary check avoids fitting. Let
\[
 \mathcal E_{g,n}=A_{g,n}-1-\sum_{j=1}^4 d_n^j g^{-j}.
\]
The expected remainder is \(O_n(g^{-5})\). The following values show
the behavior of \(g^5\mathcal E_{g,n}\):
\[
\begin{array}{c|rrr}
 &g=40&g=60&g=85\\ \hline
n=2&-0.5748239805&-0.5573173047&-0.5476248646\\
n=3&-0.4898574375&-0.4746763865&-0.4662847314\\
n=4&-0.4871451412&-0.4747020831&-0.4678525406
\end{array}
\]
The scaled residuals remain of order one and exhibit convergence,
consistent with the claimed four-term expansion.

\Addresses

\begin{thebibliography}{99}

\bibitem{AN}
A. Alexandrov and P. Norbury, 
\emph{Super volumes and KdV tau functions}, arXiv:2412.17272.

\bibitem{arbarello1996combinatorial} E. Arbarello and M. Cornalba, {\em Combinatorial and algebro-geometric cohomology classes on moduli spaces of curves}, J. Algebraic Geometry {\bf 5} (1996), 705--709.

\bibitem{CGG}
N.~K.~Chidambaram, E.~Garcia-Failde, and A.~Giacchetto,
\emph{Relations on \(\overline{\mathcal M}_{g,n}\) and the negative
\(r\)-spin Witten conjecture},
Invent. Math. \textbf{241} (2025), no.~3, 929--997.


\bibitem{NP1} N. Do and P. Norbury, {\em Topological recursion on the Bessel curve}, Comm. Num. Theor. Phys. {\bf 12} (2018): 53-73.

\bibitem{GNYZ}
J.~Guo, P.~Norbury, D.~Yang, and D.~Zagier,
\emph{Combinatorics and large genus asymptotics of the
Br\'ezin--Gross--Witten numbers},
arXiv:2412.20388v2 (2025).


\bibitem{griguolo2024asymptotics} L. Griguolo, J. Papalini, L. Russo and D. Seminara, {\em Asymptotics of Weil-Petersson volumes and two-dimensional quantum gravities}, SciPost Phys. 17, 156 (2024).

\bibitem{HLX} X. Huang, K. Liu and H. Xu, {\em Higher Weil-Petersson volumes of the moduli space of super Riemann surfaces}, arXiv:2501.05806.

\bibitem{Huang} X. Huang, {\em Large genus asymptotics of super Weil-Petersson volumes}, arXiv:2501.07848.

\bibitem{HPZ} Y. Huang, and R. Penner, and A. Zeitlin, {\em Super McShane identity}, J. Diff. Geom. {\bf 125} (2023), 509--551.

\bibitem{kimura2020jt} Y. Kimura, {\em JT gravity and the asymptotic Weil-Petersson volume}, Phys. Lett. B {\bf 811} (2020), 135989.


\bibitem{mirzakhani2013growth} M. Mirzakhani, {\em Growth of Weil-Petersson volumes and random hyperbolic surface of large genus}, J. Differential Geom. {\bf 94} (2013), 267--300.

\bibitem{MZ}
M.~Mirzakhani and P.~Zograf,
\emph{Towards large genus asymptotics of intersection numbers on
moduli spaces of curves},
Geom. Funct. Anal. \textbf{25} (2015), 1258--1289.
\href{https://arxiv.org/abs/1112.1151}{arXiv:1112.1151}.

\bibitem{norbury2023new} P. Norbury, {\em A new cohomology class on the moduli space of curves}, Geom. Top. {\bf 27} (2023), 2695--2761.

\bibitem{norbury2020enumerative} P. Norbury, {\em Enumerative geometry via the moduli space of super Riemann surfaces} J. Geom. Phys. {\bf 222} (2026), 105750.

\bibitem{stanford2020jt} D. Stanford and E. Witten, {\em JT gravity and the ensembles of random matrix theory}, Adv. Theor. Math. Phys. {\bf 24} (2020), 1475--1680.



\end{thebibliography}
\end{document}